\documentclass[11pt,a4paper,reqno]{amsart}
\usepackage[latin1]{inputenc} 
\usepackage{amsthm,amsmath,amssymb} 
\usepackage{enumerate} 
\usepackage{graphicx}

\theoremstyle{plain}
\newtheorem{theo}{Theorem}[section]

\newtheorem{lemma}[theo]{Lemma}

\theoremstyle{definition}
\newtheorem{defn}[theo]{Definition}

\newtheorem{setup}[theo]{Setup}
\newtheorem{con}[theo]{Construction}

\newcommand{\mc}[1]{\mathcal{#1}}
\newcommand{\mb}[1]{\mathbb{#1}}
\newcommand{\nib}[1]{\noindent {\bf #1}}

\newcommand{\bsize}[1]{\left| #1 \right|}

\newcommand{\bangle}[1]{\left\langle #1 \right\rangle}

\newcommand{\sub}{\subseteq}

\newcommand{\sm}{\setminus}

\newcommand{\es}{\emptyset}
\newcommand{\pl}{\partial}

\newcommand{\gG}{\gamma}

\newcommand{\lL}{\lambda}
\newcommand{\tT}{\theta}
\newcommand{\sS}{\sigma}
\newcommand{\oO}{\omega}

\newcommand{\GG}{\Gamma}
\newcommand{\DD}{\Delta}
\newcommand{\OO}{\Omega}
\newcommand{\TT}{\Theta}
\newcommand{\Ss}{\Sigma}

\def\qed{\hfill $\Box$}

\date{}
\title{Counting designs}
\author{Peter Keevash} 
\thanks{ 
Mathematical Institute, University of Oxford, Oxford, UK. 
Email keevash@maths.ox.ac.uk.
Research supported in part by ERC grant 239696. 
}

\begin{document}

\vspace*{-0.8cm}
\begin{abstract}
We give estimates on the number of combinatorial designs, which prove (and generalise)
a conjecture of Wilson from 1974 on the number of Steiner Triple Systems. 
This paper also serves as an expository treatment of our recently developed method 
of Randomised Algebraic Construction: we give a simpler proof of a special case of our result
on clique decompositions of hypergraphs, namely triangle decompositions of quasirandom graphs. 
\end{abstract}
\maketitle
\vspace*{-0.6cm}

\section{Introduction}

When does a graph $G$ have a triangle decomposition? 
(By this we mean a partition of its edge set into triangles.)
There are two obvious necessary `divisibility conditions':
the number of edges must be divisible by three,
and the degree of any vertex must be even.
We say that $G$ is \emph{tridivisible} if it satisfies 
these divisibility conditions. In 1847 Kirkman proved
that any tridivisible complete graph has a triangle decomposition;
equivalently, there is a Steiner Triple System on $n$ vertices
if $n$ is $1$ or $3$ mod $6$.
In \cite{K} we showed more generally that a tridivisible graph 
has a triangle decomposition if we assume a certain pseudorandomness condition.
In fact, we proved a more general result on clique decompositions of simplicial complexes, 
which in particular proved the Existence Conjecture for combinatorial designs.

One purpose of the current paper is to illustrate the new technique 
(Randomised Algebraic Construction) of \cite{K} in the simplified setting
of triangle decompositions; we will also prove a conjecture 
of Wilson \cite{W6} on the number of Steiner Triple Systems.
These results are proved in the next three sections,
roughly following the method of \cite{K}, 
but introducing some novelties in technique
that lead to considerable simplifications
in the case of triangle decompositions;
the material here closely follows a lecture series that the author
recently gave at the Israel Institute for Advanced Studies.
In Section \ref{rtr} we sketch an argument of Bennett and Bohman \cite{BB}
on the random greedy matching process and adapt the calculations to the version 
needed in this paper. We generalise from Steiner Triple Systems to designs in Section \ref{des}.
We conclude by noting that it remains open to obtain an asymptotic formula
for the number of designs, or even just for the number of regular graphs.

\section{Triangle decompositions}

We start by stating our result that tridivisible
pseudorandom graphs have triangle decompositions.
The pseudorandomness condition is as follows. 
Let $G$ be a graph on $n$ vertices. 
The \emph{density} of $G$ is $d(G) = |G|/\tbinom{n}{2}$. We say that $G$ is
\emph{$c$-typical} if every vertex has $(1 \pm c)d(G) n$ neighbours
and every pair of vertices have $(1 \pm c) d(G)^2 n$ common neighbours.
(We write $b \pm c$ for any real between $b - c$ and $b + c$.)

\begin{theo} \label{qrtri} 
There exists $0<c_0<1$ and $n_0 \in \mb{N}$ so that if $n \ge n_0$ 
and $G$ is a $c$-typical tridivisible graph on $n$ vertices
with $d(G) > n^{-10^{-7}}$ and $c < c_0 d(G)^{10^6}$
then $G$ has a triangle decomposition.
\end{theo}

Note that in Theorem \ref{qrtri} we allow the density to decay polynomially with $n$;
this will be important for the application in the next subsection, but in many cases of 
interest one can consider $d(G)$ and $c$ to be fixed constants independent of $n$.
One such consequence of Theorem \ref{qrtri} noted in \cite{K}
is that the standard random graph model $G(n,1/2)$ with high probability
(whp) has a partial triangle decomposition that covers all but $(1+o(1))n/4$ edges.
Indeed, deleting a perfect matching on the set of vertices of odd degree
and then at most two $4$-cycles gives a graph satisfying the hypotheses of the theorem.
This is asymptotically best possible, as whp there are $(1+o(1))n/2$
vertices of odd degree, and any set of edge-disjoint triangles must leave 
at least one edge uncovered at each vertex of odd degree.

We remark that our definition of typicality here is weaker than that used in \cite{K}.
In fact, for most of the paper we will assume the stronger version, then explain
at the end how the proof can be modified to work with the current definition.
We also make the (well-known) remark that typicality implies the standard regularity
property (for appropriate constants) that appears in Szemer\'edi's Regularity Lemma, 
but the converse is not true, as regularity allows individual vertices to behave badly,
even to be isolated.

\subsection{The number of Steiner Triple Systems}

Another purpose of our paper is to prove the following conjecture 
of Wilson \cite{W6} on the number of Steiner Triple Systems on $n$ vertices, 
i.e.\ triangle decompositions of the complete graph $K_n$; denote this by $STS(n)$.

\begin{theo} \label{wilson-conj}
If $n$ is $1$ or $3$ mod $6$, then $STS(n) = (n/e^2 + o(n))^{n^2/6}$.
\end{theo}

Note that $K_n$ is tridivisible if and only if $n$ is $1$ or $3$ mod $6$,
so $STS(n)=0$ for all other $n$. The upper bound in Theorem \ref{wilson-conj} 
was recently proved by Linial and Luria \cite{LL}, who showed that
$STS(n) \le (n/e^2 + O(\sqrt{n}))^{n^2/6}$. Our lower bound will be
$STS(n) \ge (n/e^2 + O(n^{1-a}))^{n^2/6}$ for some small $a>0$.
 
Theorem \ref{wilson-conj} will follow quite easily from Theorem \ref{qrtri}
and the semirandom method (nibble). It will be most convenient for us to apply the 
results of Bohman, Frieze and Lubetzky \cite{BFL} on the triangle removal process
(although we could make do with a simpler nibble argument,
or the argument of Bennett and Bohman \cite{BB} sketched in Section \ref{rtr}).
We say that an event $E$ holds \emph{with high probability} (whp) 
if $\mb{P}(E) = 1-e^{-\Omega(n^c)}$ for some $c>0$ as $n \to \infty$; 
note that when $n$ is sufficiently large,
by union bounds we can assume that any specified polynomial number of such events all occur. 

In the triangle removal process, we start with the complete graph $K_n$,
and at each step we delete the edges of a uniformly random triangle in the current graph.
It is shown in \cite{BFL} that whp the process persists until only $O(n^{3/2+o(1)})$ edges remain,
but we will stop at $n^{2-10^{-7}}$ edges (i.e.\ at the nearest multiple of $3$ 
to this number) so that we can apply Theorem \ref{qrtri}.
We need the following additional facts from \cite{BFL} about this stopped process:
whp the final graph is $n^{-1/3}$-typical, and when $pn^2/2$ edges remain
the number of choices for the deleted triangle is $(1 \pm n^{-2/3}) (pn)^3/6$.

\medskip

\nib{Proof of Theorem \ref{wilson-conj}.}
Consider the following procedure for constructing a Steiner Triple System on 
$n$ vertices: run the triangle removal process until $n^{2-10^{-7}}$ edges remain,
then apply Theorem \ref{qrtri} (if its hypotheses are satisfied,
which occurs in $1-o(1)$ proportion of all instances of the process).
Writing $m$ for the number of steps and $p(i) = 1-6i/n^2$,
the logarithm of the number of choices in this procedure is 
\[ L_1 = \sum_{i=1}^m (\log (p(i)^3n^3/6) \pm 2n^{-2/3})
= (n^2/6)( \log(n^3/6) - 3 \pm n^{-10^{-8}} ), \]
since $\sum_{i=1}^m \log p(i) = (1 + O(n^{-10^{-7}}\log n)) (n^2/6) \int_0^1 \log p \ dp $
and $\int_0^1 \log p \ dp = -1$. Also, for any fixed Steiner Triple System,
the logarithm of the number of times it is counted by this procedure is at most
\[L_2 =  \sum_{i=1}^m \log(p(i)n^2/6) = (n^2/6)(\log(n^2/6) - 1 \pm n^{-10^{-8}} ).\]
Therefore $\log(STS(n)) \ge L_1-L_2 = (n^2/6)(\log(n) - 2 \pm 2n^{-10^{-8}})$,
which implies the stated bound on $STS(n)$. \qed

\subsection{Strategy}

The strategy of the proof of Theorem \ref{qrtri} is encapsulated by the following setup
(we give motivation and discussion below).
We say that $J \sub G$ is \emph{$c$-bounded} if $|J(v)| < c|V(G)|$ for every $v \in V(G)$,
where $J(v) = \{u \in V(G): uv \in J\}$ is the \emph{neighbourhood} of $v$ in $J$.

\begin{setup} \label{strategy}
Suppose we have $G^* \sub G$ with a `template' triangle decomposition $T$ such that
\begin{description}
\item[Nibble] $G \sm G^*$ contains a set $N$ of edge-disjoint triangles 
 with `leave' $L := (G \sm G^*) \sm \cup N$ that is $c_1$-bounded,
\item[Cover] For any $L \sub G \sm G^*$ that is $c_1$-bounded, there is 
 a set $M^c$  of edge-disjoint triangles such that $L = (G \sm G^*) \cap (\cup M^c)$ 
 and the `spill' $S := G^* \cap (\cup M^c)$ is $c_2$-bounded,
\item[Hole] For any tridivisible $S \sub G^*$ that is $c_2$-bounded,
 there are `outer' and `inner' sets $M^o,M^i$ of edge-disjoint triangles in $G^*$ such that 
 $\cup M^o$ is $c_3$-bounded and $(S,\cup M^i)$ is a partition of $\cup M^o$,
\item[Completion] Given $L$, $M^c$, $M^o$ and $M^i$ as above, there are sets 
 $M_1$, $M_2$, $M_3$, $M_4$ of edge-disjoint triangles in $G^*$ such that $(L,\cup M_2)$ 
 is a partition of $\cup M_1$, $\cup M_3 = \cup M_4$, $M_3 \sub T$ and $M_2 \sub M_4$.
\end{description}
\end{setup}

The key step is choosing $T$ (which determines $G^*$).
We will use our method of Randomised Algebraic Construction, 
which takes a particularly simple form for triangle decompositions.
To motivate the construction, suppose that $V(G)$ is an abelian group,
and consider the set $\Ss$ of triples $xyz$ such that $x+y+z=0$.
We note that $\Ss$ is a good `model' for a triangle decomposition,
as for any $xy$ there is a unique $z$ such that $x+y+z=0$.
However, we cannot simply take $\Ss$, as not all such $xyz$ are triangles of $G$;
moreover, $x,y,z$ may not even be pairwise distinct. 

The idea of the construction is that a suitable random subset of $\Ss$ can act
as a template, which covers a constant fraction of $G$. Next we find an approximate
decomposition of the rest of $G$ by random greedy algorithms: this is accomplished by
steps \textbf{Nibble} and \textbf{Cover} of Setup \ref{strategy}.
After these steps, every edge of $G$ has been covered once or twice,
and the spill $S$ is the set of edges that have been covered twice.
Finally, we use local modifications built into the template to turn the
approximate decomposition into an exact decomposition: this is accomplished by 
steps \textbf{Hole} and \textbf{Completion} of Setup \ref{strategy}.

To motivate \textbf{Completion}, we imagine first that we have \textbf{Hole} 
and also $M^o \sub T$. Then we could delete $M^o$ and take $M^i$ 
instead, thus reducing by one the multiplicity of every edge in $S$, so that
we have a triangle decomposition of $G$. However, specifying a triangle of $T$
is very restrictive, as there are only order($n^2$) such triangles out of a total
of order($n^3$) triangles in $G$. If we had chosen $T$ uniformly at random it
would be hopeless to obtain any useful configuration formed by triangles of $T$. 
However, the algebraic structure implies that certain
configurations of triangles are dense within a sparse configuration space
(described by linear constraints). This forms the basis of a modification
procedure that replaces $M^c$, $M^o$ and $M^i$ by other sets of triangles
with the same properties, where $M_1$ plays the role of $M^c \cup M^i$,
$M_2$ of $M^o$, and each triangle $f$ of $M_2$ can be embedded in a 
small subgraph that has one triangle decomposition (part of $M_4$) using $f$
and another triangle decomposition (part of $M_3$) contained in $T$.

It is not hard to see that $G$ contains a triangle decomposition in Setup \ref{strategy}.
Indeed, we start by taking the sets $N$ provided by \textbf{Nibble}
and then the sets $M^c$ and $S$ provided by \textbf{Cover}. 
Now we note that $S = \cup T + \cup N + \cup M^c - G$ is tridivisible,
as any integer linear combination of tridivisible graphs is tridivisible.
So we can apply \textbf{Hole} to obtain $M^o$ and $M^i$. 
Then we can apply \textbf{Completion} to obtain $M_1,M_2,M_3,M_4$.
Finally, $M = N \cup M_1 \cup (M_4 \sm M_2) \cup (T \sm M_3)$
is a triangle decomposition of $G$. Thus the remainder of the proof
will be to show that we can achieve Setup \ref{strategy}.

\subsection{Template} \label{sec:template}

We choose the template as follows.

\begin{con} \label{rac}
Let $a \in \mb{N}$ be such that $2^{a-2} < |V(G)| \le 2^{a-1}$.
Let $\pi:V(G) \to \mb{F}_{2^a} \sm \{0\}$ be a uniformly random injection.
Let \[T = \{xyz \in K_3(G): \pi(x)+\pi(y)+\pi(z)=0\}
\ \text{ and } \ G^* = \cup T.\]
\end{con}

To avoid cumbersome notation, we use $xyz$ to denote either the vertex set
$\{x,y,z\}$ or the edge set $\{xy,xz,yz\}$ of a triangle.
The context determines which interpretation is intended,
e.g.\ in Construction \ref{rac} the graph $G^*$ is the (disjoint) union
of the edge-sets of the triangles in $T$.

In this subsection we will show that whp the pair $(G,G^*)$ is `typical'
(in a precise sense defined below); this will allow us to implement the
approximate decomposition in steps \textbf{Nibble} and \textbf{Cover}.
Moreover, we will show in Section \ref{sec:completion} that $G^*$ is
`linearly typical' (roughly speaking: we can count subgraph extensions 
with linear constraints on the vertices); this will imply the existence
of the local modifications used in steps \textbf{Hole} and \textbf{Completion}.

We start with some notation and preliminary observations.
Throughout we write $n=|V(G)|$. We identify $G$ with its edge set $E(G)$, 
so that $|G|$ denotes the number of edges of $G$ (rather than the number of vertices, 
as is used by some authors). We write $[n]=\{1,\dots,n\}$.
We define \[\gG = 2^{-a} n,\] and note that $1/4 < \gG < 1/2$.
We observe that if $x,y,z \in \mb{F}_{2^a} \sm \{0\}$ and $x+y+z=0$
then $x,y,z$ are pairwise distinct. We note that $+1=-1$ in $\mb{F}_{2^a}$,
so we can use $+$ and $-$ interchangeably in $\mb{F}_{2^a}$-arithmetic.
We consider $\mb{F}_{2^a}$ as a vector space over $\mb{F}_2$, 
and observe that any two nonzero elements span a subspace of dimension two.

Next we introduce the stronger typicality assumption used in \cite{K}.
We say that $G$ is \emph{$(c,h)$-typical} if 
\[|\cap_{x \in S} G(x)| = (1 \pm |S|c) d(G)^{|S|} n
\ \text{ for any } S \sub V(G)
 \text{ with }  |S| \le h.\] 
Note that being $c$-typical is essentially 
the same as being $(c,2)$-typical (up to a factor of $2$ in $c$).
For most of the paper we will assume that $G$ is $(c,16)$-typical;
at the end we will explain how the proof can be modified
to work with the weaker assumption that $G$ is $c$-typical. 

Now we define the typicality condition for $(G,G^*)$ and show that it holds whp.
Let $G^*$ be a subgraph of $G$. We say that $(G,G^*)$ is \emph{$(c,h)$-typical} if 
\[| \bigcap_{x \in S^*} G^*(x) \cap \bigcap_{x \in S \sm S^*} G(x) | 
= (1 \pm |S|c) d(G^*)^{|S^*|} d(G)^{|S|-|S^*|}n\] 
for any $S^* \sub S \sub V(G)$ with $|S| \le h$.

\begin{lemma} \label{template-typ}
whp $d(G^*) = (1 \pm 3c)\gG d(G)^3$ and $(G,G^*)$ is $(6c,16)$-typical.
\end{lemma}

The proof uses the following consequence of Azuma's inequality. 

\begin{defn} \label{def:lip2}
Let $S_n$ be the symmetric group, $f:S_n \to \mb{R}$ and $b \ge 0$.
We say that $f$ is \emph{$b$-Lipschitz} if for any 
$\sS,\sS' \in S_n$ such that $\sS = \tau \circ \sS'$ 
for some transposition $\tau \in S_n$
we have $|f(\sS)-f(\sS')| \le b$. 
\end{defn}

\begin{lemma} \label{lip2} (see e.g.\ \cite{McD})
Suppose $f:S_n \to \mb{R}$ is $b$-Lipschitz,
$\sS \in S_n$ is uniformly random and $X=f(\sS)$.
Then \[\mb{P}(|X-\mb{E}X|>t) \le 2e^{-t^2/2nb^2}.\]
\end{lemma}

\nib{Proof of Lemma \ref{template-typ}.}
We start by estimating $\mb{E}|G^*| = \sum_{e \in G} \mb{P}(e \in G^*)$.
For any $e=xy$, given $\pi(x)$ and $\pi(y)$, we have $e \in G^*$
if and only if $\pi(z) = \pi(x) + \pi(y)$ for some $z$ such that $xyz \in K_3(G)$.
Since $G$ is $(c,16)$-typical, there are $(1 \pm 2c) d(G)^2 n$ choices for $z$.
Each satisfies $\pi(z) = \pi(x) + \pi(y)$ with probability $(2^a-3)^{-1}$,
so $\mb{E}|G^*| = |G| (1 \pm 2c) d(G)^2 n (2^a-3)^{-1}$.
We can view $\pi$ as $\sS \circ \pi_0$, where $\pi_0:V(G) \to \mb{F}_{2^a} \sm \{0\}$
is any fixed injection and $\sS$ is a random permutation of $\mb{F}_{2^a} \sm \{0\}$.
Any transposition of $\sS$ affects $|G^*|$ by $O(n)$,
so by Lemma \ref{lip2} whp $d(G^*) = (1 \pm 2.1c) \gG d(G)^3$.

Similarly, we consider any $S^* \sub S \sub V(G)$ with $|S| \le 16$,
write $Y = \bigcap_{x \in S^*} G^*(x) \cap \bigcap_{x \in S \sm S^*} G(x)$,
and estimate $\mb{E}|Y| = \sum_{y \in V(G)} \mb{P}(y \in Y)$.
For any $y \in \cap_{x \in S} G(x)$, given $\pi(y)$ and $\pi(x)$ for all $x \in S$,
we have $y \in Y$ if and only if for all $x \in S^*$
there is $xyz_x \in K_3(G)$ such that $\pi(z_x) = \pi(x) + \pi(y)$.
Since $G$ is $(c,16)$-typical, there are $(1 \pm |S|c) d(G)^{|S|} n$ choices for $y$.
By excluding $O(1)$ choices of $y$ we can assume 
$\pi(x) + \pi(y) \ne \pi(x')$ for all $x,x' \in S$.
Then for each $x \in S^*$ there are $(1 \pm 2c) d(G)^2 n$ choices for $z_x$,
and for any set of choices, with probability $(1+O(1/n)) 2^{-a|S^*|}$
they all satisfy $\pi(z_x) = \pi(x) + \pi(y)$.
This gives \[\mb{E}|Y| = O(1) + (1 \pm |S|c) d(G)^{|S|} n 
\cdot ((1 \pm 2c) d(G)^2 n)^{|S^*|} \cdot (1+O(1/n)) 2^{-a|S^*|}.\]
Any transposition of $\sS$ affects $|Y|$ by $O(1)$,
so by Lemma \ref{lip2} whp $|Y| = (1 \pm (3|S|+1)c) d(G)^{|S|} (\gG d(G)^2)^{|S^*|} n
= (1 \pm 6|S|c) d(G^*)^{|S|} d(G)^{|S|-|S^*|} n$. \qed

\medskip

Since $d(G^*) = (1 \pm 3c)\gG d(G)^3$ and $1/4 < \gG < 1/2$
we have $0.24d(G)^3 < d(G^*) < 0.51 d(G)$ for small $c$.
Also, as $(G,G^*)$ is $(6c,16)$-typical
we can deduce that $G \sm G^*$ is $50c$-typical. 
Indeed, for any $v \in V(G)$ we have 
\begin{align*} 
& |(G \sm G^*)(v)| = (1 \pm c)d(G)n - (1 \pm 6c)d(G^*)n \\
& = (d(G)-d(G^*))n \pm 6c(d(G)+d(G^*))n 
= (1 \pm 20c)d(G \sm G^*)n. 
\end{align*}
Furthermore, for any $u,v \in V(G)$ we estimate 
$|(G \sm G^*)(u) \cap (G \sm G^*)(v)|$ as
\begin{align*} 
&  |G(u) \cap G(v)| - |G^*(u) \cap G(v)| - |G(u) \cap G^*(v)| + |G^*(u) \cap G^*(v)| \\
& = (1 \pm 2c)d(G)^2 n - 2(1 \pm 12c)d(G)d(G^*) n + (1 \pm 12c)d(G^*)^2 n \\
& = (d(G)-d(G^*))^2 n \pm 12c(d(G)+d(G^*))^2 n
= (1 \pm 50c) d(G \sm G^*)^2 n. \end{align*}
Applying the following theorem, we deduce \textbf{Nibble} with $c_1 = (50c)^{1/4}$.

\begin{theo} \label{nibble-deg}
There are $b_0>0$ and $n_0 \in \mb{N}$ so that if $n>n_0$, $n^{-0.1} < b < b_0$ 
and $G$ is a $b$-typical graph on $n$ vertices with $d(G)>b$, then there is 
a set $N$ of edge-disjoint triangles in $G$ such that $L = G \sm \cup N$ is $b^{1/4}$-bounded.
\end{theo}

We remark that the parameters in Theorem \ref{nibble-deg} are not very sharp:
we have just fixed some convenient values that suffice for our purposes.
Similar results are well-known, but we are not aware of any reference
that implies the theorem as stated, so we will sketch a proof in Section \ref{rtr}.

For convenient reference, we give here the values of some 
other parameters that will be used below: 
\begin{gather*} c_2 = 10^2 c_1 d(G)^{-6}, \quad c_3 = 10^{20} c_2 d(G)^{-50} \\
c_4 = 10^{20} c_3 d(G)^{-100} \quad \text{ and } c_5 = 10^{10} c_4 d(G)^{-180}.
\end{gather*} 
The tightest constraint on $c$ that will be required in our calculations is
$100 c_5 = 10^{54} (50c)^{1/4} d(G)^{-336} < 10^{-6} d(G)^{180}$; 
this holds for small $c_0$ if $c < c_0 d(G)^{3000}$.
(This is the bound we need if $G$ is $(c,16)$-typical,
but if $G$ is $c$-typical we need the stronger bound in Theorem \ref{qrtri}.)

\subsection{Cover}

Consider the following random greedy algorithm. 
Let $L = \{e_i: i \in [t]\}$ (with edges ordered arbitrarily).
Let $M^c = \{T_i: i \in [t]\}$ be triangles such that $T_i$ consists of $e_i$
and two edges of $G^*$, and is chosen uniformly at random from all such
triangles that are edge-disjoint from all previous choices;
if there is no available choice for $T_i$ then the algorithm aborts.

To analyse the algorithm we require a concentration inequality.
We say that a random variable $Y$ is \emph{$(\mu,C)$-dominated},
if there are constants $\mu_1,\dots,\mu_m$ with $\sum_{i=1}^m \mu_i < \mu$,
and we can write $Y = \sum_{i=1}^m Y_i$, such that $|Y_i| \le C$ for all $i$,
and conditional on any given values of $Y_j$ for $j<i$ we have $\mb{E}|Y_i|<\mu_i$.
The following lemma follows easily from Freedman's inequality \cite{F}
(see \cite[Lemma 2.7]{K}). 

\begin{lemma} \label{dom} 
If $Y$ is $(\mu,C)$-dominated then $\mb{P}(|Y|>2\mu) < 2e^{-\mu/6C}$.
\end{lemma}

Sometimes we will use a modified inequality with $2$ replaced by $1+c$.
We also note that if the $Y_i$ are independent (not necessarily identically distributed)
indicator variables we recover a version of the Chernoff bound for (pseudo)binomial
variables (where better concentration is known).
For the following lemma, we recall that $L$ is $c_1$-bounded, 
where $c_1 = (50c)^{1/4}$, and that $c_2 = 10^2 c_1 d(G)^{-6}$.

\begin{lemma} \label{cover-leave} 
whp the algorithm to choose $M^c$ does not abort,
and $S := G^* \cap (\cup M^c)$ is $c_2$-bounded.
\end{lemma}

\nib{Proof.}
For $i \in [t]$ we let $\mc{B}_i$ be the bad event that 
$S_i := G^* \cap (\cup_{j<i} T_j)$ is not $c_2$-bounded.
We define a stopping time $\tau$ be the smallest $i$ for which $\mc{B}_i$ holds
or the algorithm aborts, or $\infty$ if there is no such $i$. 
It suffices to show whp $\tau=\infty$. 

We fix $t_0 \in [t]$ and bound $\mb{P}(\tau=t_0)$ as follows.
For any $i<t_0$, since $\mc{B}_i$ does not hold, $S_i$ is $c_2$-bounded.
Writing $e_i = v_iv'_i$, we can bound the number of excluded choices for $T_i$ 
by $c_2 n < |G^*(v_i) \cap G^*(v'_i)|/2$, 
so at most one half of the triangles on $e_i$ are excluded.

Next we fix $e = vv' \in G^*$, and estimate 
$r_e := \sum_{i \le t_0} \mb{P}'(e \sub T_i)$,
where $\mb{P}'$ denotes the conditional probability
given the choices made before step $i$.
We compare $r_e$ to the expected number of times
that $e$ would be covered if we chose all triangles independently.
To be precise, we let \[E_e := \sum_{i \le t_0} \mb{P}(e \sub T'_i),\]
where each $T'_i$ is a uniform random triangle consisting of $e_i$
and two edges of $G^*$, and $(T'_i: i \in [t])$ are independent.
By the bound on excluded choices, 
$\mb{P}'(e \sub T_i) < 2 \mb{P}(e \sub T'_i)$,
so $r_e < 2E_e$.

The $i$th summand in $E_e$ is only nonzero when $e_i \cap e \ne \es$.
As $L$ is $c_1$-bounded, the number of such $i$
is at most $|L(v)|+|L(v')|<2c_1n$.
Also, for each $i$ such that $e_i \cup e$ spans a triangle, we have 
\[\mb{P}(e \sub T'_i) = |G^*(v_i) \cap G^*(v'_i)|^{-1} < 2 d(G^*)^{-2} n^{-1}.\]
Therefore $E_e < 4c_1 d(G^*)^{-2} < c_2/4$.

Finally, fix $v \in V(G)$ and consider $X = |S_{t_0}(v)| = \sum_{i \le t_0} X_i$,
where $X_i = \sum_{v \in e \in G^*} 1_{e \sub T_i}$. We have $|X_i| \le 2$ and
\[\sum_{i \le t_0} \mb{E}'(X_i) = \sum_{i \le t_0} \sum_{v \in e \in G^*} 
\mb{P}'(e \sub T_i) = \sum_{v \in e \in G^*} r_e < c_2 n/2.\]
By Lemma \ref{dom} we have $\mb{P}(X \ge c_2 n) < 2e^{-c_2 n/24}$. 
Taking a union bound over $i \le t_0 \le t$, whp $|S(v)| < c_2 n$,
i.e.\ $S$ is $c_2$-bounded and $\tau=\infty$. \qed

\medskip

Below we will require several more random greedy algorithms similar to that above.
One could formulate an abstract general lemma to cover all cases (see \cite[Lemma 4.11]{K}), 
but here we will prefer the more intuitive approach of identifying the key principles of the proof,
so that it will be clear how it may be adapted to future instances.
For a general random greedy algorithm, we identify some desired boundedness conclusion,
then at each step of the algorithm, assuming that boundedness has not failed, we show that 
at most one half (say) of the choices of the required configuration have been excluded. 
Then for each edge $e$ in the underlying graph $H$ we estimate the expected number $E_e$ 
of times that $e$ would be covered if we chose all configurations independently.
If $E_e<b/4$ and the configurations have constant size (not depending on $n$)
then the graph of all covered edges is whp $b$-bounded.

We record some estimates that are useful for such arguments.
Suppose $H$ is a small fixed graph ($|H| \le 500$ say),
$F \sub V(H)$ and $\phi$ is an embedding of $H[F]$ in $G^*$.
We call $E=(\phi,F,H)$ an \emph{extension}.
Let $X_E(G^*)$ be the number of embeddings $\phi^*$ of $H$ in $G^*$ that restrict to $\phi$ on $F$.
We suppose that $E$ is \emph{$16$-degenerate}, meaning that we can construct the embedding one
vertex at a time, so that at each step we add a vertex adjacent to at most $16$ existing vertices.
As $(G,G^*)$ is $(6c,16)$-typical, when we add a vertex adjacent to $t \le 16$ existing vertices, 
there are $(1 \pm 6tc)d(G^*)^t n$ choices. Multiplying these estimates, we obtain 
the following estimate for $X_E(G^*)$.

\begin{lemma} \label{ext}
Suppose $E=(\phi,F,H)$ is a $16$-degenerate extension with $|H| \le 500$. 
Then \[X_E(G^*) = (1 \pm 7|H|c) d(G^*)^{|H \sm H[F]|} n^{|V(H)|-|F|}.\] 
\end{lemma}

Now suppose that we wish to exclude embeddings $\phi^*$ that use some edge in $J$,
which is $c$-bounded. Fix $e \in H \sm H[F]$ and consider the embeddings $\phi^*$
with $\phi^*(e) \in J$. If $e \cap F \ne \es$ there are at most $cn$ choices for the
embedding of $e$ then at most $n^{|V(H)|-|F|-1}$ choices for the remainder of $\phi^*$.
If $e \cap F = \es$ there are at most $cn^2$ choices for the embedding of $e$
then at most $n^{|V(H)|-|F|-2}$ choices for the remainder of $\phi^*$.
Thus at most $|H|cn^{|V(H)|-|F|}$ choices of $\phi^*$ are excluded,
which is a negligible fraction of $X_E(G^*)$.

\section{Integral relaxations}

In this section we establish \textbf{Hole}.
Our first step is to consider an integral relaxation,
in the following sense. Instead of thinking of
$(S,\cup M^i)$ as a partition of $\cup M^o$, we think of $S$
as a weighted sum of edge sets of triangles, where triangles in $M^o$
have weight $1$ and triangles in $M^i$ have weight $-1$. 
We can express this by the equation $\Phi A = S$,
where $\Phi$ is the corresponding $\pm 1$-vector indexed by triangles, 
and $A$ is the inclusion matrix of triangles against edges,
i.e.\ $A_{fe}=1_{e \sub f}$ for any edge $e$ and triangle $f$.
It is straightforward to show that this equation has a solution if
we allow $\Phi$ to have any integer weights on triangles
(see \cite{GJ,W4,W5} for more general results).

It will be more convenient to work with linear maps rather than matrices.
For any graph $H$ we define $\mb{Z}$-linear boundary/shadow maps 
$\pl_j : \mb{Z}^{K_i(H)} \to \mb{Z}^{K_j(H)}$ for $i \ge j \ge 0$
by $\pl_j(e) = \sum \tbinom{e}{j}$ for $e \in K_i(H)$,
i.e.\ for $J \in \mb{Z}^{K_i(H)}$ and $f \in K_j(H)$
we define $\pl_j(J)_f = \sum_{f \sub e \in K_i(H)} J_e$.
For example, if $J \in \mb{Z}^H$ then $\pl_1(J) \in \mb{Z}^{V(H)}$
is defined by $\pl_1(J)_v = \sum_{v \in e \in H} J_e$.

It will also be notationally convenient to identify vectors with (generalised) sets.
It is standard to identify $v \in \{0,1\}^X$ with the set $\{x \in X: v_x=1\}$.
Similarly, we can identify $v \in \mb{N}^X$ with the multiset in $X$ in which $x$ 
has multiplicity $v_x$ (for our purposes $0 \in \mb{N}$). We also apply similar notation 
and terminology as for multisets to vectors $v \in \mb{Z}^X$ (`intsets').
Here our convention is that `for each $x \in v$' means that $x$ is considered $|v_x|$ times 
in any statement or algorithm, and has a sign attached to it (the same as that of $v_x$);
we also refer to $x$ as a `signed element' of $v$. For $v \in \mb{Z}^X$ we write $v=v^+-v^-$, 
where $v^+_x = \max\{v_x,0\}$ and $v^-_x = \max\{-v_x,0\}$ for $x \in X$. Given $J \in \mb{N}^G$
and $v \in V(G)$, we define $J(v) \in \mb{N}^{V(G)}$ by $J(v)_u = 1_{uv \in G} J_{uv}$.
Then we can extend the definition of boundedness to multigraphs:
$J$ is $c$-bounded if $|J(v)| < cn$ for every $v \in V(G)$.

With this notation, our integral relaxation of \textbf{Hole} is expressed by 
the following lemma (in which $K_n$ denotes the complete graph on $V(G)$);
for \textbf{Hole} we will need the additional properties
that $\Phi(f) = 0$ for any $f \in K_3(K_n) \sm K_3(G^*)$,
and $\Phi(f) \in \{0,1,-1\}$ for all $f \in K_3(G^*)$,
as then we can write $\Phi = M^o - M^i$.

\begin{lemma}\label{intrelax}
There is $\Phi \in \mb{Z}^{K_3(K_n)}$ with $\pl_2 \Phi = S$ 
such that $\pl_2 \Phi^+$ is $100c_2$-bounded.
\end{lemma}

\nib{Proof.}
We will construct $\Phi = \Phi_0 + \Phi_1 + \Phi_2$ 
such that $J^0 = S - \pl_2 \Phi_0$, $J^1 = J^0 - \pl_2 \Phi_1$, 
$J^2 = J^1 - \pl_2 \Phi_2$ satisfy $\pl_i J^i = 0$ for $i=0,1,2$.
Recalling that $S$ is tridivisible, each $J^i$ will be tridivisible,
in the `intgraph' sense: i.e.\ $\sum_e J^i_e$ is divisible by $3$
and $\sum_u J^i_{uv}$ is divisible by $2$ for all $v$.

\medskip

\emph{Step 0:}
For $\Phi_0$, we choose $|S|/3$ independent uniformly random triangles in $K_n$; 
then $J^0 = S - \pl_2 \Phi_0$ satisfies $\pl_0 J^0 = 0$. For each vertex $v$, 
the number of these triangles containing $v$ is binomial with mean $|S|/n < c_2 n/2$,
so by the Chernoff bound whp $\pl_2 \Phi_0$ is $1.1c_2$-bounded.

\medskip

\emph{Step 1:}
We let $J^* = \pl_1 J^0$, so $\pl_0 J^* = 2 \pl_0 J^0 = 0$, i.e. $|J^{*+}|=|J^{*-}|$. Note 
for all $x \in V(G)$ that $J^*_x$ is even, as $J^0$ is tridivisible, and $|J^*_x| < 1.1c_2 n$. 
We fix an arbitrary sequence $((x^+_i,x^-_i) : i \in [|J^{*+}|/2])$ so that each $x \in V(G)$ 
occurs $J^{*+}_x/2$ times as some $x^+_i$ and $J^{*-}_x/2$ times as some $x^-_i$. 
For each $i$ we choose $a_ib_i \sub V(G) \sm \{x^+_i,x^-_i\}$ independently uniformly at random,
and let $\Phi_1 = \sum_{i \in [|J^{*+}|/2]} (\{x^+_ia_ib_i\}-\{x^-_ia_ib_i\})$;
then  $J^1 = J^0 - \pl_2 \Phi_1$ satisfies $\pl_1 J^1 = 0$.

We claim that whp $\pl_2 \Phi_1^\pm$ are $8c_2$-bounded.
To see this, we first fix any $e \in K_n$ and estimate the expected contributions to $e$ 
from each step $i$, according to whether $e$ contains $x^+_i$, $x^-_i$, or neither.
Each endpoint of $e$ occurs at most $0.6c_2 n$ times as $x^\pm_i$, and for such $i$
we cover $e$ with probability $2/(n-2)$, so the expected contribution to $(\pl_2 \Phi_1^\pm)_e$ 
from all such $i$ is at most $2.5c_2$. At any other step, we cover $e$ with probability 
$\tbinom{n-2}{2}^{-1}$, so the total expected contribution to $(\pl_2 \Phi_1^\pm)_e$ from these
steps is at most $1.1c_2$. Now, for each vertex $v$, summing over its incident edges,
$|\pl_2 \Phi_1^\pm(v)|$ are both $(4c_2 n,1)$-dominated, so the claim holds by Lemma \ref{dom}.

\medskip

\emph{Step 2:}
We start by fixing an arbitrary expression $J^1 = \sum_{C \in \mc{C}_0} C$, where each $C$ 
is a closed walk in $G^*$ with edge weights alternating between $1$ and $-1$, and there are no 
cancellations, i.e.\ every edge appears in the sum only with weight $1$ or only with weight $-1$. 
As is well-known, such an expression may be found by a greedy algorithm: each $C$ can be obtained
by following an arbitrary alternating walk on the signed elements of $J^1$ until we return 
to our starting point using an edge with the opposite sign to that of the first edge, 
whereupon we add $-C$ to $J^1$ and repeat the procedure.
(We note that this argument leads to a convenient shortcut for triangle decompositions,
but does not generalise to hypergraph decompositions.)

\begin{figure}
\begin{center}
\includegraphics[scale=0.8]{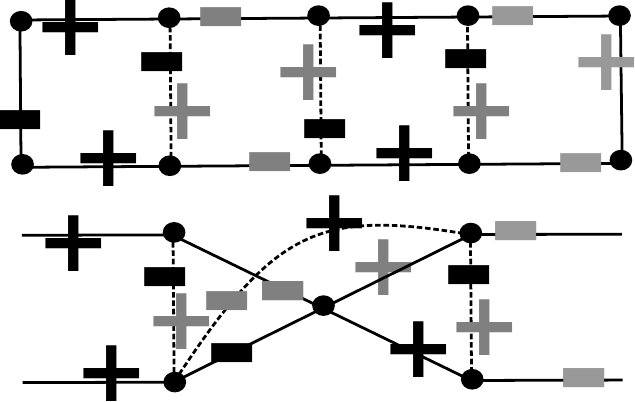}
\caption{Decomposing even signed cycles.}
\end{center}
\end{figure}

Next we express each $C \in \mc{C}_0$ as a sum of signed four-cycles in the complete graph $K_n$ 
on $V(G)$, where we write each closed walk of length $2m$ as a chain of $m-1$ signed four-cycles, 
using the identity (see Figure 1)
\begin{align*} 
& \sum_{i=1}^{m-1} (-1)^i (\{x_ix_{i+1}\} - \{x_{i+1}y_{i+1}\} + \{y_{i+1}y_i\} - \{y_ix_i\}) \\
& = \{x_1y_1\} + (-1)^m \{x_my_m\}
+ \sum_{i=1}^{m-1} (-1)^i \{x_ix_{i+1}\} + \sum_{i=1}^{m-1} (-1)^i \{y_iy_{i+1}\} .  
\end{align*}
This identity can be used as is if $x_i \ne y_i$ for $i \in [m]$.
For each $i$ such that $x_i=y_i$, we note that $1<i<m$, $x_{i-1} \ne y_{i-1}$, $x_{i+1} \ne y_{i+1}$,
and $x_{i+1} \ne y_{i-1}$, so we can replace the four-cycles for summands $i-1$ and $i$ by 
\begin{align*} 
& (-1)^{i-1}(\{x_{i-1}x_i\}-\{x_ix_{i+1}\}+\{x_{i+1}y_{i-1}\}-\{y_{i-1}x_{i-1}\}, \ \text{ and } \\
& (-1)^i(\{x_{i+1}y_{i-1}\}-\{y_{i-1}y_i\}+\{y_iy_{i+1}\}-\{y_{i+1}x_{i+1}\}).
\end{align*}
Thus we can write $J^1 = \sum_{C \in \mc{C}} C$, where each summand is a signed four-cycle in $K_n$.
Furthermore, the above construction has the property that for each $v \in V(G)$ and $w \in \{-1,1\}$
we use at most $3|J^{1+}(v)| < 24c_2 n$ edges at $v$ with weight $w$.

For each $C = \{ab\}-\{bc\}+\{cd\}-\{da\} \in \mc{C}$ we choose 
$x \in V(G) \sm \{a,b,c,d\}$ independently uniformly at random,
and add $\{xab\}-\{xbc\}+\{xcd\}-\{xda\}$ to $\Phi_2$; 
then $\pl_2 \Phi_2 = \sum_{C \in \mc{C}} C = J^1$.
Let $\GG$ denote the multigraph formed by
summing $\{xa,xb,xc,xd\}$ over all such $C$.
For any $e \in K_n$, at most $48c_2 n$ elements of $\mc{C}$
can contribute to $\GG_e$, so $\mb{E}\GG_e < 49c_2 n$.
Then for any $v$, summing over its incident edges, 
$|\GG(v)|$ is $(49c_2 n,4)$-dominated, so by Lemma \ref{dom} (modified)
whp $\GG$ is $50c_2$-bounded. Defining $\Phi = \Phi_0 + \Phi_1 + \Phi_2$,
we have $\pl_2 \Phi = S$ and $\pl_2 \Phi^+$ is $100c_2$-bounded. \qed

\medskip

\begin{figure}
\begin{center}
\includegraphics[scale=0.8]{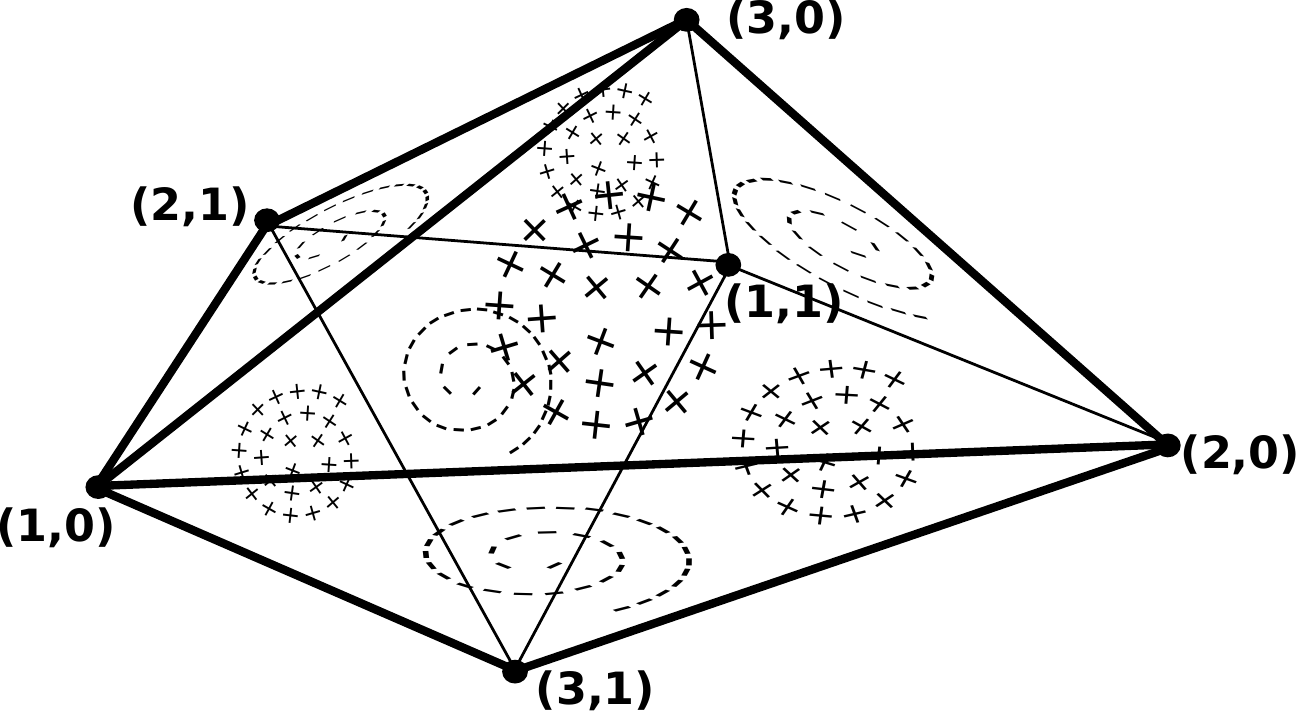}
\caption{An octahedron with signed triangles.}
\end{center}
\end{figure}

To obtain \textbf{Hole}, we will modify $\Phi$ using 
the following `octahedral' configurations (see Figure 2).
Consider a copy of $K_{2,2,2}$, the complete tripartite graph with $2$ points in each part,
with parts $\{(j,0), (j,1)\}$ for $j \in [3]$. We denote its triangles by 
$\{f_x: x \in \{0,1\}^3\}$, where $f_x = \{ (j,x_j): j \in [3]\}$. 
The sign of $f_x$ is $s(f_x) = (-1)^{\sum x}$. Thus each edge is in one triangle of each sign.
Defining $\OO = \sum_{x \{0,1\}^3} s(f_x) \{f_x\} \in \mb{Z}^{K_3(K_{2,2,2})}$, we see that $\pl_2 \OO = 0$. 
This gives a method to eliminate any signed triangle $f$ from $\Phi$ without altering $\pl_2 \Phi$:
we add some copy of $\OO$ with the opposite sign to $f$ in which (say) $f_{000} = f$, 
thus replacing $f$ by seven other signed triangles that have the same total $2$-shadow. 
Similarly (and more importantly), we can eliminate any pair of triangles $f,f'$
that have opposite sign and share an edge $e$, replacing $f,f'$ by six other
signed triangles that have the same total $2$-shadow and do not use $e$.
We apply this method in the following two-phase algorithm.

\medskip

\emph{Octahedral Elimination Algorithm (Phase I).}
We eliminate all triangles in $\Phi$, according to a random greedy algorithm,
where in each step we consider some original signed element $f$ of $\Phi$,
and choose an octahedral configuration $\OO_f$ to replace $f$.
We refer to edges of $\OO_f$ not in $f$ as \emph{new} edges, and choose $\OO_f$
uniformly at random subject to the new edges belonging to $G^*$ and being disjoint 
from $\pl_2 \Phi^+$ and all new edges from previous steps.

\medskip

Let $\Phi'$ denote the result of Phase I (if it does not abort).
Then $\pl_2 \Phi' = \pl_2 \Phi = S$, and we can write $\pl_2 \Phi'{}^+ = \pl_2 \Phi^+ + \GG$,
where $\GG$ is the graph of new edges, and every signed element of $\Phi'$
contains at most one edge of $\pl_2 \Phi^+$.

\medskip

\emph{Octahedral Elimination Algorithm (Phase II).}
We replace all signed edges apart from those in $S$ and $\GG$. To do this, 
we fix a sequence $\mc{S}$ of pairs of signed elements of $\Phi'$, so that 
(i) for each $ff' \in \mc{S}$, there is some $e \in \pl_2 \Phi^+$ such that 
$f$ and $f'$ both contain $e$, and $f$ and $f'$ have opposite signs, and
(ii) the multiset consisting of all $e$ as in (i) is $\pl_2 \Phi^-$.
Now we eliminate each $ff' \in \mc{S}$, according to a random greedy algorithm,
by subtracting some copy $\OO_{ff'}$ of $\OO$ with $f_{000} = f$ and $f_{001} = f'$,
or vice versa, depending on the signs.
We refer to edges of $\OO_{ff'}$ not in $f$ or $f'$ as \emph{new} edges, and choose $\OO_{ff'}$
uniformly at random subject to the new edges belonging to $G^*$ and being distinct 
from $\pl_2 \Phi^+ \cup \GG$ and all new edges from previous steps.

\medskip

Let $\Psi$ denote the result of this algorithm (if it does not abort)
and $\GG'$ the graph of new edges for Phase II.
Then $\pl_2 \Psi = S$ and $\pl_2 \Psi^- = \GG \cup \GG' \sub G^*$.
This implies $\Psi(f) = 0$ for any $f \in K_3(K_n) \sm K_3(G^*)$,
and $\Psi(f) \in \{0,1,-1\}$ for all $f \in K_3(G^*)$, so $\Psi = M^o - M^i$, 
where $M^o$ and $M^i$ are as in \textbf{Hole}, once we have verified the boundedness condition.

\begin{lemma}
whp the Octahedral Elimination Algorithm produces 
$M^o$ and $M^i$ as in \textbf{Hole}.
\end{lemma}

\nib{Proof.}
We first show that whp $\GG$ is $c'_2$-bounded, where $c'_2 = 10^5 c_2 d(G^*)^{-9}$. 
The proof follows the discussion after the proof of Lemma \ref{cover-leave},
where a configuration for $f$ consists of the new edges of some $\OO_f$. 
By Lemma \ref{ext}, at each step, the number of choices of $\OO_f$ with all new edges 
belonging to $G^*$ (with no excluded configurations) is $(1 \pm 60c) d(G^*)^9 n^3$.
Assuming that the graph of previous new edges is $c'_2$-bounded,
as $\pl_2 \Phi^+$ is $100 c_2$-bounded, the number of excluded  configurations 
is at most $10c'_2 n^3$, which is less than half of the total. 
Next, for each $e \in G^*$, we consider separately the contributions to $E_e$,
according to whether $e$ intersects $f$ in $0$ or $1$ vertex
(there is no contribution to new edges from triangles containing $e$).
There are at most $600 c_2 n$ signed elements of $\Phi$ that intersect $e$ in $1$ vertex.
For each of these, a random configuration covers $e$ with probability 
at most $3n^2 / (1 - 60c) d(G^*)^9 n^3$, so the total contribution to $E_e$
from such elements is at most $2000 c_2 d(G^*)^{-9}$.
Also, $\Phi$ has at most $100 c_2 n^2$ signed elements,
and for each one that is disjoint from $e$ the contribution to $E_e$
is at most $6n / (1 - 60c) d(G^*)^9 n^3$, so the total contribution
from such elements is at most $1000 c_2 d(G^*)^{-9}$.
We obtain $E_e < 3000 c_2 d(G^*)^{-9}$, which implies the claimed bound on $\GG$. 

Next we claim that whp $\GG'$ is $c''_2$-bounded, where $c''_2 = 20c'_2 d(G^*)^{-7}$. 
The argument is very similar to that given for $\GG$.
Now a configuration for $ff'$ consists of the new edges of some $\OO_{ff'}$. 
By Lemma \ref{ext}, at each step, the number of choices of $\OO_{ff'}$ with all new edges 
belonging to $G^*$ (with no excluded configurations) is $(1 \pm 50c) d(G^*)^7 n^2$.
Assuming that the graph of previous new edges is $c''_2$-bounded,
as $\pl_2 \Phi^+ \cup \GG$ is $2c'_2$-bounded,
the number of excluded configurations is at most $10c''_2 n^2$, 
which is less than half of the total. Next, for each $e \in G^*$, 
we consider separately the contributions to $E_e$
according to whether $e$ intersects $f \cup f'$ in $0$ or $1$ vertex
(there is no contribution to new edges if $e \sub f \cup f'$).

First we consider those $ff' \in \mc{S}$ that intersect $e$ in $1$ vertex $x$.
There are two choices for $x \in e$. If $x \in f \cap f'$ then there are
at most $200 c_2 n$ choices for $f \cap f' \in \pl_2 \Phi^+ \cup \pl_2 \Phi^-$, 
which determines $f$ and $f'$. If $\{x\} = f \sm f'$
then there are at most $|\GG(x)| < c'_2 n$ choices for $f$, and so $f'$.
The same bound applies if $\{x\} = f' \sm f$,
so there are at most $5c'_2 n$ such $ff'$.
Each contributes at most $2n / (1 - 50c) d(G^*)^7 n^2$ to $E_e$, so the total 
contribution from such $ff'$ is at most $11 c'_2 d(G^*)^{-7}$.
Also, $|\mc{S}| = |\pl_2 \Phi^-| < 100 c_2 n^2$, and for each $ff' \in \mc{S}$
with $e \cap (f \cup f') = \es$ the contribution to $E_e$ is at most 
$2 / (1 - 50c) d(G^*)^7 n^2$, so the total contribution
from such elements is at most $300 c_2 d(G^*)^{-7}$.
We obtain $E_e < 12 c'_2 d(G^*)^{-7}$, which implies the claimed bound on $\GG'$.
Recalling that $d(G^*) > 0.24d(G)^3$ and $c_3 = 10^{20} c_2 d(G)^{-50}$
we see that $\cup M^o = \pl_2 \Psi^+ = S \cup \GG \cup \GG'$ is $c_3$-bounded,
so we have the required properties for \textbf{Hole}. \qed

\section{Completion} \label{sec:completion}

For \textbf{Completion}, we divide the analysis into two parts.
Firstly, we will determine what conditions on $M_1$ and $M_2$
enable us to find $M_3$ and $M_4$. 
Secondly, we will show that the sets $M^c$, $M^o$ and $M^i$ from \textbf{Cover} and \textbf{Hole}
can be modified to give $M_1$ and $M_2$ satisfying the required conditions.
For convenient notation we suppress the embedding $\pi:V(G) \to \mb{F}_{2^a}$ whenever 
we do not need to refer to it, instead thinking of $V(G)$ as a subset of $\mb{F}_{2^a}$.

\subsection{Shuffles.}
Suppose we have a set $M_2$ of edge-disjoint triangles in $G^*$, and we want 
to find sets $M_3$ and $M_4$ of edge-disjoint triangles in $G^*$
such that $\cup M_3 = \cup M_4$, $M_3 \sub T$ and $M_2 \sub M_4$.
Our basic building blocks (`shuffles') will be edge-disjoint subgraphs of $G^*$,
each having two different triangle decompositions, one only using triangles in $T$, 
and the other including any specified triangle of $M_2$. Then the unions over all
blocks of the two triangle decompositions will give $M_3$ and $M_4$ as required.

We define the shuffles as follows. 
Fix $x=(x_1,x_2,x_3) \in \mb{F}_{2^a}^3$ and $t=(t_1,t_2) \in \mb{F}_{2^a}^2$
such that $\{x_1,x_2,x_3,t_1,t_2\}$ is linearly independent over $\mb{F}_2$. 
Let $\bangle{x}$ be the subspace of $\mb{F}_{2^a}$ generated by $\{x_1,x_2,x_3\}$. 
The $xt$-shuffle $S_{xt}$ is the complete tripartite graph with parts 
$t_i + \bangle{x} = \{t_i+y: y \in \bangle{x}\}$, $i \in [3]$, where $t_3:=t_1+t_2$.
If $S_{xt} \sub G^*$ then it has a triangle decomposition $M_{3xt}$ only using triangles in $T$:
take all triangles $y_1y_2y_3$ where each $y_i \in t_i+\bangle{x}$ and $y_1+y_2+y_3=0$.
We define another triangle decomposition $M_{4xt}$ of $S_{xt}$ by translating each triangle
of $M_{3xt}$ by $(x_1,x_2,x_3)$, i.e.\ $M_{4xt}$ consists of all triangles $y_1y_2y_3$ 
where each $y_i \in t_i+\bangle{x}$ and $x_1+x_2+x_3+y_1+y_2+y_3=0$.

To construct $M_3$ and $M_4$, we choose shuffles according to a random greedy algorithm,
where in each step we consider some $z_1z_2z_3 \in M_2$, and choose some shuffle 
$S_{xt} \sub G^*$ such that $z_i = t_i+x_i$ for all $i \in [3]$.
We will see in Lemma \ref{shuffle-count} that the Randomised Algebraic Construction 
is whp such that there are many choices for such a shuffle.
This is the most important property of the construction, and it would not hold 
if we had chosen the template to be a uniformly random set of edge-disjoint triangles;
in fact the expected number of shuffles (or any `shuffle-like' configuration) would be $o(1)$.
First we identify a property that we need for triangles in $M_2$ so that the required shuffles 
exist and can be chosen to be edge-disjoint. We say that $z_1z_2z_3$ is \emph{octahedral} 
if $z_1+z_2+z_3 \ne 0$ and there is a copy $K'$ of $K_{2,2,2}$ in $G$ such that 
$\pi(K')$ has parts $\{z_1,z_2+z_3\}$, $\{z_2,z_1+z_3\}$ and $\{z_3,z_1+z_2\}$; 
we call $K'$ the \emph{associated octahedron} of $z_1z_2z_3$. We assume
\begin{itemize}
\item[(P1)] all triangles in $M_2$ are octahedral, 
with edge-disjoint associated octahedra.
\end{itemize}

\begin{lemma}\label{shuffle-count}
Under the random choice of $\pi$ used in the definition of $T$,
whp for any octahedral $z_1z_2z_3$ there are $(1 \pm 200c) d(G)^{180} \gG^{18} 2^{2a}$
shuffles $S_{xt} \sub G^*$ such that $t_i+x_i=z_i$ for $i \in [3]$.
\end{lemma}

\nib{Proof.} 
We can write the number of such shuffles as a sum of indicator variables
$X=\sum 1_{E(K,\ell,x,t)}$, where the sum ranges over all $(K,\ell,x,t)$ 
such that $K$ is a copy of $K_{8,8,8}$ in $G$ containing 
the associated octahedron $K'$ of $z_1z_2z_3$, 
$\ell$ is a bijective labelling of each part of $K$ by $\mb{F}_2^3$,
we let $E(K,\ell,x,t)$ be the event that $\pi(w)=t_i+\ell(w)\cdot x$ 
for all $i \in [3]$ and $w$ in the $i$th part of $K$,
and we assume $\ell$ is consistent with $K'$, in that
$\ell(\pi^{-1}(z_i))=e_i$ and $\ell(\pi^{-1}(z_i+z_j))=e_i+e_j$ for $\{i,j\} \sub [3]$.

As $G$ is $(c,16)$-typical, there are $(1 \pm 181c) d(G)^{180} n^{18}$ choices of $(K,\ell)$. 
There are $2^{2a} - O(n)$ choices of $t$, which determines $x$ given $z$,
as only $O(n)$ choices of $t$ are excluded by the condition that
$\{x_1,x_2,x_3,t_1,t_2\}$ is linearly independent over $\mb{F}_2$:
there are $O(1)$ possible linear relations between them, 
and each such relation is linearly independent or contradictory 
to the system $t_i+x_i=z_i$ for $i \in [3]$ (as $z_1+z_2+z_3 \ne 0$),
so is satisfied by at most $2^a$ choices of $t$.
Given $(K,\ell,x,t)$, conditional on $\pi\vert_{K'}$,
we have $\mb{P}(E(K,\ell,x,t)) = (1+O(1/n)) 2^{-18a}$.
Therefore $\mb{E}X = (1 \pm 182c) d(G)^{180} \gG^{18} 2^{2a}$.

Also, any transposition $\tau$ of $\pi$ affects $X$ by at most $100 \cdot 2^a$.
To see this, we estimate the number of shuffles containing $z_1z_2z_3$ and any fixed 
$v \in \mb{F}_{2^a} \sm \{z_1,z_2,z_3,z_1+z_2,z_1+z_3,z_2+z_3\}$. Consider any $j \in [3]$,
$b \in \mb{F}_2^3 \sm \{e_j, (1,1,1)-e_j\}$, and the equations 
$t_i + b \cdot x = v$ and $t_i+x_i=z_i$ for $i \in [3]$ in $(t,x)$.
We have four linearly independent constraints, so there are at most $2^a$ solutions.
Including multiplicative factors for $i$, $b$ and $\tau$ gives the required bound. 
Now by Lemma \ref{lip2} whp $X = (1 \pm 200c) d(G)^{180} \gG^{18} 2^{2a}$. \qed

\subsection{Linear extensions.}
We digress to note a more general estimate for future reference.
Suppose $H$ is a graph, $y=(y_i:i \in [g])$ are variables, and for all $v \in V(H)$ we have distinct
linear forms $L_v(y) = c_v + \sum_{i \in S_v} y_i$ for some $c_v \in \mb{F}_{2^a}$ and $S_v \sub [g]$. 
We call $E=(L,H)$ a \emph{linear extension} with \emph{base} $F = \{v \in V(H): S_v=\es\}$. 
Let $X_E(G^*)$ be the number of \emph{$L$-embeddings} of $H$, i.e.\ embeddings $\phi$ of $H$ in $G^*$ 
such that for some $y \in \mb{F}_{2^a}^g$ we have $\phi(v) = L_v(y)$ for all $v \in V(H)$.
The above argument (see also \cite[Lemma 5.15]{K})
gives the following formula analogous to that obtained for shuffles.

\begin{lemma} \label{linext}
Let $E=(L,H)$ be a $16$-degenerate linear extension with $|H| \le 500$. Suppose
\begin{itemize}
\item $H$ has a triangle decomposition $M$ such that for each $xyz \in M$ we have $L_x+L_y=L_z$,
\item The incidence matrix of $\{S_v: v \in V(H)\}$ has full column rank $g \ge 1$.
\end{itemize}
Then \[X_E(G^*) =  (1 \pm 1.1|H|c) d(G)^{|H \sm H[F]|} \gG^{|V(H) \sm F|} 2^{ga}.\] 
\end{lemma}

\subsection{Shuffle algorithm.}
Recalling our general framework for random greedy algorithms, we want to show 
that, of the potential shuffles $S_{xt}$ with $t_i+x_i=z_i$ for $i \in [3]$,
at most half are excluded due to sharing an edge with a previous shuffle, assuming 
some boundedness condition on the graph $\GG$ of new edges from previous shuffles. 
We classify the potential restrictions according to the label of the shuffle edge
involved, which is specified by some $\{j,k\} \sub [3]$ and $b_j,b_k \in \mb{F}_2^3$ 
such that $b_j \notin \{(e_j,(1,1,1)-e_j)$ or $b_k \notin \{(e_k,(1,1,1)-e_k)$ 
(here we do not consider edges of the associated octahedra: 
these are already determined, and edge-disjoint by (P1).)
For any $v_jv_k \in G^*$, the shuffles excluded because of mapping the given
labelled shuffle edge to $v_j v_k$ are given by the $(x,t)$-solutions of 
the system $\mc{S}$ of equations $t_j + b_j \cdot x = v_j$, $t_k + b_k \cdot x = v_k$ 
and $t_i+x_i=z_i$ for $i \in [3]$. There may be $0$, $1$ or $2^a$ solutions.
We can ignore the case of $0$ solutions, as it does not exclude anything.
For the cases with $1$ solution, we can bound the number of excluded choices 
by the number of edges covered by all shuffles, which is $192|M_2|$.

It remains to consider the case that $\mc{S}$ has $2^a$ solutions, which occurs when 
one of the equations is redundant, due to being a linear combination of the other equations. 
There are a constant number of linear combinations, and each constrains $(v_j,v_k)$
to lie on a line, as may be seen from general considerations of linear algebra,
or simply by enumerating the possibilities: wlog $t_k + b_k \cdot x = v_k$ is redundant,
due to 
\begin{enumerate}[(i)]
\item $b_k=e_k$ and $v_k=z_k$, 
\item $b_k=(1,1,1)-e_k$ and $v_k=z_1+z_2+z_3-z_k$,
\item $b_j+b_k=e_j+e_k$ and $v_j+v_k=z_j+z_k$, 
\item $b_j+b_k=e_i$ and $v_j+v_k=z_i$, where $[3]=\{i,j,k\}$.
\end{enumerate}
In cases (i) and (ii) where $v_k$ is fixed, assuming that $\GG$ is $c_5$-bounded,
there are at most $c_5 n$ choices for $v_j$ such that $v_jv_k \in \GG$. 
In cases (iii) and (iv) we need an additional boundedness condition: 

We say that $\GG$ is \emph{linearly $c_5$-bounded} if $\GG$ is $c_5$-bounded 
and also contains at most $c_5 2^a$ edges from any line
of the form $\{(x_1+\mu,x_2+\mu): \mu \in \mb{F}_{2^a}\}$.

We also need similar conditions so that we can avoid the associated octahedra;
writing $\DD$ for the union of all associated octahedra of triangles in $M_2$,
we will ensure that
\begin{itemize}
\item[(P2)] $\DD$ is linearly $c_4$-bounded. 
\end{itemize}
Then the total number of excluded shuffles is at most
$192(|M_2| + (c_4+c_5) 2^{2a}) < 200c_5 2^{2a}$, which is less than half of the total.

Next we fix $e \in G^*$ and estimate $E_e$. To do so, we fix $b_j,b_k$ as above,
write $e=v_jv_k$ and estimate the sum over $z_1z_2z_3 \in M_2$
of the probability $p$ that a random shuffle $S_{xt}$ with $t_i+x_i=z_i$ for $i \in [3]$ 
satisfies $t_j + b_j \cdot x = v_j$ and $t_k + b_k \cdot x = v_k$. For fixed $z_1z_2z_3$,
if the system $\mc{S}$ as above has $N$ solutions then 
$p = N/(1 \pm 200c) d(G)^{180} \gG^{18} 2^{2a}$. When $N=1$ the total contribution 
is at most $|M_2|/(1 - 200c) d(G)^{180} \gG^{18} 2^{2a} < 1.1c_4 d(G)^{-180} \gG^{-18}$.
If $N=2^a$ then $(z_1,z_2,z_3)$ is constrained to 
lie in a certain plane (this can be seen by linear algebra, or by considering
each possibiity as above: e.g.\ in case (iii) the plane is $v_j+v_k=z_j+z_k$).
Thus we see the final property that we need from $M_2$: 
\begin{itemize}
\item[(P3)] $M_2$ contains at most $c_4 2^a$ elements $z_1z_2z_3$ from any 
\emph{basic plane} of the form $b \cdot z = v$ where $b \in \mb{F}_2^3 \sm \{0\}$. 
\end{itemize}
(Note that by (P1) we can assume $v \ne 0$ in (P3).)
Then the total contribution is at most $c_4 2^a \cdot 2^a /(1 - 200c) d(G)^{180} \gG^{18} 2^{2a}$.
Summing over $\{j,k\}$, $b_j$ and $b_k$, we can estimate $E_e < 250 c_4 d(G)^{-180} \gG^{-18} = c_5/4$.
Applying Lemma \ref{dom} as in the proof of Lemma \ref{cover-leave}, we deduce that whp the boundedness
assumptions on $\GG$ used above do not fail (linear boundedness follows in the same
way as boundedness), and so the algorithm does not abort. This completes the analysis
of the first part of \textbf{Completion}: given $M_1$ and $M_2$ as in \textbf{Completion},
under the conditions (P1--P3) on $M_2$, we can find $M_3$ and $M_4$ as in \textbf{Completion}.

\subsection{Octahedral Elimination Algorithm.}
To complete the proof of \textbf{Completion}, and so of the theorems, it remains to show that 
we can find $M_1$ and $M_2$ satisfying the conditions (P1--P3). We apply a similar two-phase
algorithm to that used in \textbf{Hole}.

\medskip

\emph{Phase I.} We start with $\Phi = M^c + M^i - M^o$, 
so $\pl_2 \Phi = L$, $\pl_2 \Phi^+ = \cup(M^c \cup M^i)$, $\pl_2 \Phi^- = \cup M^o$.
Next we eliminate all triangles in $\Phi$ according to a random greedy algorithm,
where in each step we consider some original signed element $f$ of $\Phi$,
and choose an octahedral configuration $\OO_f$ to replace $f$.
We say that a triangle $f'$ of $\OO_f$ is \emph{far} if $|f' \cap f| \le 1$,
and that $\OO_f$ is \emph{valid} if 
(i) none of its triangles are template triangles,
with the possible exception of $f$, and
(ii) all of its far triangles are octahedral, 
and their associated octahedra share edges only in $\OO_f$, 
in which case we denote their union by the extended configuration $\OO^+_f$. 
We say that an edge of $\OO^+_f$ not in $f$ is \emph{new}, 
and choose a valid $\OO_f$ uniformly at random subject to the new edges being 
distinct from all new edges from previous steps. 

\medskip

Let $\Phi'$ denote the result of Phase I (if it does not abort).
We have $\pl_2 \Phi' = \pl_2 \Phi = L$, and
writing $\GG$ for the graph of new edges,
every signed element of $\Phi'$ is either a far triangle
consisting of three edges of $\GG$,
or is not far and consists of two edges of $\GG$
and one edge of $\pl_2 \Phi^+$.

\medskip

\emph{Phase II.}
Now we will eliminate all triangles of $\Phi'$ apart from those that contain an edge of $L$
or were far in the previous modification procedure. We partition all such triangles into 
a sequence $\mc{S}$ of pairs of signed elements of $\Phi'$, so that 
for each $ff' \in \mc{S}$, there is some $e \in \pl_2 \Phi^+$ such that 
$f$ and $f'$ both contain $e$, and $f$ and $f'$ have opposite signs.
We eliminate each $ff' \in \mc{S}$, according to a random greedy algorithm,
by subtracting some copy $\OO_{ff'}$ of $\OO$ with $f_{000} = f$ and $f_{001} = f'$,
or vice versa, depending on the signs. Now we say that $\OO_{ff'}$ is \emph{valid}
if all of its triangles apart from $f$ and $f'$ are octahedral,
and their associated octahedra share edges only in $\OO_{ff'}$, 
in which case we denote their union by the extended configuration $\OO^+_{ff'}$. 
We refer to edges of $\OO^+_{ff'}$ not in $f$ or $f'$ as \emph{new} edges, 
and choose a valid $\OO_{ff'}$ uniformly at random subject to the new edges being 
distinct from $\GG$ and all new edges from previous steps.

\medskip

Let $\Psi$ denote the result of this algorithm (if it does not abort)
and $\GG'$ the graph of new edges for Phase II.
Since $\pl_2 \Psi = \pl_2 \Phi = L$, defining $M_1=\Psi^+$ and $M_2 = \Psi^-$, 
we see that $\cup M_2 = \GG \cup \GG'$ and $\cup M_1 = L \cup \GG \cup \GG'$,
so $(L,\cup M_2)$ is a partition of $\cup M_1$. The following lemma
completes the proof of \textbf{Completion}, and so of the theorems,
under the assumption that $G$ is $(c,16)$-typical.

\begin{lemma}
whp $M_2$ satisfies (P1), (P2) and (P3).
\end{lemma}

\nib{Proof.}
To analyse Phase I, we first estimate the number of choices for
an extended configuration on a triangle $f$. This can be described by 
the linear extension $(\OO^+_f,L)$, where $\OO^+_f$ is as above, we have 
variables $z=(z_1,z_2,z_3)$, which we also use to label the vertices of $\OO_f \sm f$, 
we define $L_x = x$ for all $x \in \OO_f$, and define $L_x$ for all other $x$ as 
required for the far triangles in $\OO_f$ to be octahedral, i.e.\ 
in the associated octahedron for a triangle $abc$,
the linear forms on the two vertices in each of the three parts are
$\{L_a,L_b+L_c\}$, $\{L_b,L_c+L_a\}$ and $\{L_c,L_a+L_b\}$.
By Lemma \ref{linext} whp $G^*$ is such that for any triangle $f$ in $\Phi$, 
there are $(1 \pm 60c) d(G)^{45} \gG^{15} 2^{3a}$ valid choices of $\OO_f$.
Here we also use the fact that for any triangle $abc$ of $\OO_f$ other 
than $f$ there are only $2^{2a}$ solutions to $L_a(z)+L_b(z)+L_c(z)=0$.
The precise exponents of $d(G)$ and $\gG$ (which are not important for 
the argument) may be easily calculated from the observation that adding 
an octahedron to a triangle adds $3$ new vertices and $9$ new edges, 
and $\OO^+_f$ is the composition of $5$ such extensions.

Next we claim that whp the graph $\GG$ of new edges is 
linearly $c'_3$-bounded, where $c'_3 = 400 c_3 d(G)^{-45} \gG^{-15}$.
We assume this bound on the current graph of new edges and estimate how many 
configurations are excluded. Consider any edge $uu'$ of the extended configuration.
Suppose first that $uu' \cap f = \es$. If $L_u(y)+L_{u'}(y)$ is not constant, 
then for any $vv' \in G^*$ the number of $L$-embeddings with $L_u(y)=v$ and $L_{u'}(y)=v'$ 
is at most $2^a$. There are at most $45(|M^c|+|M^i|+|M^o|) < 100c_3 n^2$ choices for a previous 
new edge $vv'$, so this excludes at most $100c_3 n^2 2^a$ configurations. 
On the other hand, if $L_u(y)+L_{u'}(y)$ is constant, then $L_u(y)$ and $L_{u'}(y)$
are constrained to lie on a basic line; there are at most $c'_3 2^a$ choices for $vv'$ 
by linear boundedness, and each such $vv'$ excludes at most $2^{2a}$ configurations. 
The latter estimate also applies to the case when one of $u$ or $u'$ is in $f$.
Summing these bounds over all $uu'$, we see that 
fewer than half of the total configurations are excluded.

Next we fix any edge $uu'$ of the extended configuration, any $vv' \in G^*$,
and estimate the sum over $f \in \Phi$ of the probability $p$ that a random
configuration satisfies $L_u(y)=v$ and $L_{u'}(y)=v'$. If $uu' \cap f = \es$ 
and $L_u(y)+L_{u'}(y)$ is not constant, 
then $p< 2^a / (1 - 60c) d(G)^{45} \gG^{15} 2^{3a}$ for any $f \in \Phi$.
There are at most $c_3 n^2$ choices for $f$, 
so the total contribution is at most $2c_3 d(G)^{-45} \gG^{-15}$.
Otherwise, if $L_u(y)+L_{u'}(y)$ is constant or one of $u$ or $u'$ is in $f$,
then one vertex of $f$ is specified by $(L_u(y),L_{u'}(y))$. 
For example, writing $f=abc$, in the associated octahedron for $az_2z_3$,
if $u=z_2$ and $u'=a+z_2$ then $a$ is specified by $(L_u(y),L_{u'}(y))$.
Then there are at most $2c_3 n$ choices for $f$ (as $\cup M^o$ is $c_3$-bounded).
For each such $f$ we have a contribution of at most 
$2^{2a} / (1 - 60c) d(G)^{45} \gG^{15} 2^{3a}$, so again 
the total contribution is at most $2c_3 d(G)^{-45} \gG^{-15}$.
Summing these bounds over all $uu'$ we can estimate
$E_{vv'} < 100 c_3 d(G)^{-45} \gG^{-15} = c'_3/4$.
Applying Lemma \ref{dom} as in the proof of Lemma \ref{cover-leave}, 
we deduce the claimed bound on $\GG$.

We also claim that whp there are at most $2c'_3 2^a$ far triangles 
in any basic plane $\Pi = \{z: b \cdot z = v\}$. 
To see this, we first consider the contribution from the template triangles
$\Pi^* = \Pi \cap T$. Since $z_1+z_2+z_3=0$ is linearly independent 
or contradictory to the defining equation of $\Pi$ we have $|\Pi^*| \le 2^a$.
Summing $E_{vv'} < c'_3/4$ over an edge $vv'$ in each triangle of $\Pi^*$,
by Lemma \ref{dom} whp $\Pi$ contains at most $c'_3 2^a$ template triangles.
Now fix any far non-template triangle $f'$ of the extended configuration,
any $g \in K_3(G^*)$, and estimate the sum over $f \in \Phi$ of the probability $p$
that a random configuration satisfies $L_{f'}(y)=g$. If $f' \cap f = \es$ 
then as $f'$ is non-template it determines the configuration, 
so $p < 1 / (1 - 60c) d(G)^{45} \gG^{15} 2^{3a}$,
giving a total contribution of at most $2c_3 d(G)^{-45} \gG^{-15} n^{-1}$.
Otherwise, $f'$ determines one of the associated octahedra,
so specifies one vertex of $f$, for example, writing $f=abc$,
if $f' = \{z_2,a+z_2,a+z_3\}$ then $a$ is specified.
Then there are at most $2c_3 n$ choices for $f$; for each such $f$ we have 
$p < 2^a / (1 - 60c) d(G)^{45} \gG^{15} 2^{3a}$, so again
the total contribution is at most $2c_3 d(G)^{-45} \gG^{-15} n^{-1}$.
Summing over $f'$ and applying Lemma \ref{dom} as in the proof 
of Lemma \ref{cover-leave}, we deduce the claimed bound on $\Pi$. 
This completes the analysis of Phase I.

To analyse Phase II, we first estimate the number of choices for
an extended configuration on a pair $ff'$. This can be described by 
the linear extension $(\OO^+_{ff'},L)$, where $\OO^+_{ff'}$ is as above, we have 
variables $z=(z_1,z_2)$, which we also use to label the vertices of $\OO_{ff'} \sm (f \cup f')$, 
we define $L_x = x$ for all $x \in \OO_{ff'}$, and define $L_x$ for all other $x$ as 
required for the triangles in $\OO_{ff'}$ other than $f$ and $f'$ to be octahedral.
The linear forms are distinct, as $f$ and $f'$ are not template triangles.
By Lemma \ref{linext} there are $(1 \pm 60c) d(G)^{53} \gG^{20} 2^{2a}$ valid $\OO_{ff'}$.
Again, the precise exponents of $d(G)$ and $\gG$ are not important for the argument,
but are straightforward to calculate: e.g.\ $\gG$ appears with exponent $20$, 
as $\OO_{ff'}$ and each of the $7$ associated octahedra adds $3$ new vertices to the extension,
but $4$ vertices in the associated octahedra belong to the base of the extension, 
being the third vertex of the template triangle containing an edge in $f$ or $f'$ other than $e$.

We claim that whp $\GG'$ is linearly $c_4/2$-bounded. 
The argument is very similar to that given above for $\GG$.
Assuming this bound on the current graph of new edges, 
one can show similarly to before that fewer than half 
of the total configurations are excluded.
We also need to estimate the sum over $f \in \Phi$ of the probability 
that a random configuration satisfies $L_u(y)=v$ and $L_{u'}(y)=v'$,
for any $uu'$ in the extended configuration and $vv' \in G^*$.
For most choices for $uu'$ the required bound follows as before, 
but there is an additional case, 
namely when $uu' \cap (f \cup f') = \es$ and $L_u(y)+L_{u'}(y)$ is constant, 
it may be that no vertex of $f \cup f'$ is specified by $(L_u(y),L_{u'}(y))$,
but instead some pair (not $e$) is constrained to lie on a basic line.
For example, writing $f=abc$ and $f'=abc'$, if $u = b+z_1$ and $u'=c+z_1$
then $(L_u(y),L_{u'}(y))$ specifies $b+c$, but not $b$ or $c$.
In this case, we use the fact that $\GG$ is linearly $c'_3$-bounded
to see that there are at most $c'_3 2^a$ choices for $ff'$.
Each such $ff'$ contributes at most $2^a / (1 - 60c) d(G)^{53} \gG^{20} 2^{2a}$,
giving a total contribution of at most $2c'_3 d(G)^{-53} \gG^{-20}$.
Summing over all $uu'$ we estimate $E_{vv'} < 100c'_3 d(G)^{-53} \gG^{-20}$,
so the claim follows from Lemma \ref{dom}.

Finally, $M_2$ satisfies the conditions (P1--P3): 
indeed, (P1) holds by definition of the extended configurations 
and random greedy algorithms, (P2) holds as $\DD \sub \GG \cup \GG'$,
and (P3) holds as whp $\Psi$ has at most $c_4 n$ triangles in any basic plane:
this holds for the new triangles in this algorithm by the same argument as for $\Phi'$,
and we may include the far triangles from the previous algorithm in this estimate. \qed

\medskip

This completes the proof of \textbf{Completion}, and so of Theorem \ref{qrtri},
under the assumption that $G$ is $(c,16)$-typical. 
The following modification proves the theorem under the assumption that $G$ is $c$-typical.
It is well-known that $G$ is $c^{1/50}$-regular (say) in the `Szemer\'edi sense'
(see e.g.\ \cite[Theorem 2.2]{K2}).
For $S \sub V(G)$ with $|S| \le 16$, say that $S$ is \emph{good}
if $|\cap_{x \in S} G(x)| = (1 \pm |S|c^{1/100}) d(G)^{|S|} n$, otherwise \emph{bad}.
As $G$ is $c$-typical, there are no bad sets of size $1$ or $2$. By regularity,
for $k \le 15$, any good $k$-set is contained in at most $c^{1/100} n$ bad $(k+1)$-sets.
The proof of Lemma \ref{template-typ} still applies if we assume that all subsets of $S$ are good.
Thus we can avoid using bad sets with negligible changes to the calculations.

\section{Random triangle removal} \label{rtr}

In this section we sketch a proof of Theorem \ref{nibble-deg},
by describing how to apply the analysis of random greedy hypergraph matching 
by Bennett and Bohman \cite{BB} (we choose this for simplicity, but there
are several other alternative approaches).

Consider the random triangle removal algorithm starting with $G$ rather than $K_n$.
The intuition is that after $i$ steps the remaining graph $G^i$ should look like a
random subgraph of $G$ where edges are retained independently with probability $p=1-3t$,
where $t=i/|G|$. For any $e \in G$ let $T_e(i)$ denote the number of triangles of $G^i$
containing $e$ and $Q(i)$ denote the total number of triangles of $G^i$.
By assumption, $T_e(0) = (1 \pm b)D$ for all $e \in G$, where $D = d(G)^2 n$.
Note also that $Q(0) = (1 \pm b)|G|D/3$. 

We will show that for $0 \le i \le (1-b^{1/4})|G|/3$ whp
\[Q(i) = |G|Dp^3/3 \pm e_q \ \text{ and } \ 
T_e(i) = Dp^2 \pm e_d \text{ for all } e \in G^i,\] 
where $e_q = 2(1-3\log p)^2 b|G|D$ and $e_d = 2(1-3\log p)b^{2/3} D$.
We restrict attention to the upper bounds, as the lower bounds are similar.
A convenient reformulation is to show negativity of shifted variables 
\[Q^+(i) = Q(i) - |G|Dp^3/3 - e_q \ \text{ and } \ 
T_e^+(i)=T_e(i)-Dp^2-e_d.\] 
This follows whp from martingale concentration inequalities (e.g.\ \cite{F})
after we verify the following `trend' and `boundedness' hypotheses, 
supposing that the required estimates for $Q$ and $T_e$ hold at previous steps 
(i.e.\ $i<\tau$, where the stopping time $\tau$ is the first step 
where any of the required estimates fails, or $\infty$); 
we use primes to denote conditional expectation
given the history of the process.

\medskip

\emph{Trend hypothesis:} If $Q^+(i) \ge -b|G|D$ then $\mb{E}'Q^+(i+1) \le Q^+(i)$;
if $T_e^+(i) \ge -b^{2/3}D$ then $\mb{E}'T_e^+(i+1) \le T_e^+(i)$.

\medskip

\emph{Boundedness hypothesis:} $(Q^+(i+1)-Q^+(i))^2|G|p\log n < (b|G|D)^2$,
and $-\TT < T_e^+(i+1)-T_e^+(i) < \tT$ with $\tT<\TT/10$
and $\tT \TT |G|p\log n < (b^{2/3}D)^2$.

\medskip

We start with the boundedness hypothesis, which holds with room to spare.
Indeed, $|Q^+(i+1)-Q^+(i)|=O(n)$, so $(Q^+(i+1)-Q^+(i))^2|G|p\log n=O(n^4\log n)$,
whereas $(b|G|D)^2 > b^2 d(G)^3 n^6/4 > n^{5.5}/4$. Also, for any $e \in G^i$ we have
$-1 \le T_e(i+1)-T_e(i) \le 0$. The change in $Dp^2$ is $O(Dp/|G|)$ and in $e_d$
is $O(b^{2/3}D/p|G|)$, so we can take $\TT = 2$ and $\tT = O(p+b^{2/3}/p) D/|G|$.
Then $\tT \TT |G|p\log n = O(p^2+b^{2/3}) D\log n = O(n\log n)$,
whereas $(b^{2/3}D)^2 = b^{4/3} d(G)^4 n^2 > n^{1.4}$.

Next consider the trend hypothesis for $Q$. 
Conditional on the required estimates at step $i$, if $Q^+(i) \ge -b|G|D$
we have
\begin{align*} & \mb{E}'[Q(i+1)-Q(i)]  = - Q(i)^{-1} \sum_{e \in G^i} T_e(i)^2 + O(1) \\
& = -9Q(i)/|G|p \pm 2|G|pe_d^2/Q(i) + O(1) \\
& \le - 3Dp^2 - 9e_q/|G|p + 9bD/p +  (6+o(1))e_d^2/Dp^2.
\end{align*}
Also, as $e_q$ is increasing, the one-step change 
in $-|G|Dp^3/3-e_q$ is at most $(1+O(1/|G|))3Dp^2$. As $p \ge b^{1/4}$, 
we deduce 
\begin{align*} & \mb{E}'[Q^+(i+1)-Q^+(i)] \le - 9e_q/|G|p + 9bD/p + 7e_d^2/Dp^2 \\
& \le -18(1-3\log p)^2 bD/p + 9bD/p + 28(1-3\log p)^2 b^{4/3}D/p^2 \le 0.
\end{align*}

Now consider the trend hypothesis for $T_e$.
Conditional on the required estimates at step $i$, if $T_e^+(i) \ge -b^{2/3}D$
we have 
\begin{align*} & \mb{E}'[T_e(i+1)-T_e(i)] = - Q(i)^{-1} T_e(i) 2(Dp^2 \pm e_d  + O(1)) \\
& \le - 2Q(i)^{-1}(Dp^2+e_d-b^{2/3}D) (Dp^2-e_d-O(1)) \\
& \le - 6(Dp^2)^2/|G|Dp^3 + (6+o(1))(b^{2/3}D)(Dp^2)/|G|Dp^3 + O(e_q/|G|^2p^2).
\end{align*}
Also, the one-step change in $-Dp^2$ is at most $(1+O(1/|G|))6Dp/|G|$,
and in $-e_d$ is at most $-(1+O(1/|G|p))18b^{2/3}D/|G|p$. We deduce 
\[\mb{E}'[T_e^+(i+1)-T_e^+(i)] \le -10b^{2/3}D/|G|p +  O(\log^2 p) bD/|G|p^2 \le 0.\]

Thus the required estimates hold, i.e.\ for $0 \le i \le (1-b^{1/4})|G|/3$ 
whp $Q(i) = |G|Dp^3/3 \pm e_q$ and $T_e(i) = Dp^2 \pm e_d$ for all $e \in G^i$.

Now we apply the same method to deduce the boundedness conclusion of
Theorem \ref{nibble-deg}. For $0 \le i \le (1-b^{1/4})|G|/3$ we show whp 
\[ |G^i(v)| = p|G(v)| \pm e_v \ \text{ for any vertex } v ,\]
where $e_v = 2b^{1/3} d(G)n$.
The upper bound $G_v^+(i) := |G^i(v)| - p|G(v)| - e_v \le 0$ 
will follow whp after we verify the following two conditions.

\medskip

\emph{Trend hypothesis:} If $G_v^+(i) \ge -b^{1/3}d(G)n$ then $\mb{E}'G_v^+(i+1) \le G_v^+(i)$.

\medskip

\emph{Boundedness hypothesis:} $-\TT < G_v^+(i+1)-G_v^+(i) < \tT$ with $\tT<\TT/10$
and $\tT \TT |G|p\log n < (b^{1/3}d(G)n)^2$.

\medskip

For the boundedness hypothesis, we can take $\TT=2$ and $\tT = 3|G(v)|/|G|$,
so $\tT \TT |G|p\log n = 6p|G(v)|\log n = O(n\log n)$, 
whereas $(b^{1/3}d(G)n)^2 > b^3 n^2 \ge n^{1.7}$. 
For the trend hypothesis, if $G_v^+(i) \ge -b^{1/3}d(G)n$ we have
\begin{align*} 
& \mb{E}'[|G^{i+1}(v)|-|G^i(v)|] = - \sum_{e: v \in e \in G^i} T_e(i)/Q(i) \\
& \le - (p|G(v)|+e_v-b^{1/3}d(G)n)(Dp^2-e_d)/(|G|Dp^3/3+e_q) \\
& \le -3|G(v)|/|G| - 3(e_v-b^{1/3}d(G)n)/|G|p
+ O(e_d) \tfrac{|G(v)|}{|G|Dp^2} + O(e_q) \tfrac{|G(v)|}{|G|^2Dp^3}.
\end{align*}
The one-step change in $-p|G(v)| - e_v$ is at most $3|G(v)|/|G|$, so
\[ pn\ \mb{E}'[G_v^+(i+1)-G_v^+(i)] 
\le -6b^{1/3} + O(e_d/Dp) + O(e_q/|G|Dp^2) \le 0,\]
as $e_d/Dp = O(b^{2/3}p^{-1}\log p)$, 
$e_q/|G|Dp^2 = O(bp^{-2}\log^2 p)$ and $p \ge b^{1/4}$.

Thus, letting $N$ be the set of triangles removed during the process,
whp $L = G \sm \cup N$ is $b^{1/4}$-bounded.

\section{The number of designs} \label{des}

In this section we generalise Theorem \ref{wilson-conj} to estimate the number of designs.
We start by describing the results of \cite{K} on the existence of designs. 
Let $D$ be a $q$-graph (i.e.\ a set of subsets of size $q$) of a set $X$ of size $n$. 
We say that $D$ is a \emph{design} with parameters $(n,q,r,\lL)$ 
if every subset of $X$ of size $r$ belongs to exactly $\lL$ elements of $S$.
Note that if $q=3$, $r=2$, $\lL=1$ then $D$ is a Steiner Triple System.
The necessary divisibility conditions generalise in a straightforward way:
if $D$ exists then $\tbinom{q-i}{r-i}$ must divide $\lL \tbinom{n-i}{r-i}$ for every 
$0 \le i \le r-1$; to see this, fix any $i$-subset $I$ of $X$ and consider the sets in $D$ 
that contain $I$. In \cite{K} we proved the `Existence Conjecture', which states that these
divisibility conditions are also sufficient for the existence of a design with parameters 
$(n,q,r,\lL)$, assuming $q,r,\lL$ are fixed and $n>n_0(q,r,\lL)$ is large.
We will generalise Theorem \ref{wilson-conj} as follows.

\begin{theo} \label{count}
For any $q,r,\lL$ there is $n_0$ such that if $n>n_0$ and
$\tbinom{q-i}{r-i} \mid \lL \tbinom{n-i}{r-i}$ for all $0 \le i \le r-1$,
writing $Q=\tbinom{q}{r}$ and $N=\tbinom{n-r}{q-r}$, 
the number $D(n,q,r,\lL)$ of designs with parameters $(n,q,r,\lL)$ satisfies
\[D(n,q,r,\lL) = \lL!^{-\tbinom{n}{r}} ( (\lL/e)^{Q-1} N + o(N))^{\lL Q^{-1} \tbinom{n}{r}}.\]
\end{theo}

The proof of Theorem \ref{count} follows that of Theorem \ref{wilson-conj}:
the lower bound generalises the argument given earlier in this paper,
and the upper bound generalises that of Linial and Luria \cite{LL}.

We start with the lower bound.
In the same way as a Steiner Triple System can be viewed as a triangle decomposition of $K_n$,
we can view a design with  parameters $(n,q,r,\lL)$ as a $K^r_q$-decomposition of $\lL K^r_n$,
where $K^r_q$ denotes the complete $r$-graph on $q$ vertices and $\lL K^r_n$ denotes the 
multi(hyper)graph in $K^r_n$ in which each edge has multiplicity $\lL$.
To generalise Theorem \ref{qrtri}, we first need to define the divisibility and typicality 
conditions for general $r$-graphs (we will omit multiplicities in the
definitions, as we do not need them for our application here). 

For $S \sub V(G)$, the \emph{neighbourhood} $G(S)$ 
is the $(r-|S|)$-graph $\{f \sub V(G) \sm S: f \cup S \in G\}$.
We say that $G$ is \emph{$K^r_q$-divisible} if $\tbinom{q-i}{r-i}$ divides $|G(e)|$ 
for any $i$-set $e \sub V(G)$, for all $0 \le i \le r$. 
We say that $G$ is \emph{$(c,h)$-typical}
if there is some $p>0$ such that for any set $A$ of $(r-1)$-subsets of $V(G)$
with $|A| \le h$ we have $\bsize{\cap_{S \in A} G(S)} = (1 \pm c) p^{|A|} n$.

Now we can state the $r$-graph generalisation of Theorem \ref{qrtri}.
When $d(G)$ is at least a constant independent of $n$
this follows from \cite[Theorem 1.4]{K};
the same proof shows that $d(G)$ can decay polynomially in $n$.

\begin{theo} \label{hyp-decomp} 
For any $q>r \ge 1$ there are $c_0,a \in (0,1)$ and $h,\ell,n_0 \in \mb{N}$ so that 
if $n \ge n_0$ and $G$ is a $K^r_q$-divisible $(c,h)$-typical $r$-graph on $n$ vertices 
with $d(G)>n^{-a}$ and $c < c_0 d(G)^\ell$ then $G$ has a $K^r_q$-decomposition.
\end{theo}

In the proof of Theorem \ref{count} it is more convenient to count designs together with a choice
for each $e \in K^r_n$ of a bijection between the $\lL$ copies of $e$ in $\lL K^r_n$ and the $\lL$
sets of the design containing $e$; we will refer to such a structure as an \emph{edge-labelled design}
with parameters $(n,q,r,\lL)$ and denote their number by $D^*(n,q,r,\lL)$.
As $D^*(n,q,r,\lL) = \lL!^{\tbinom{n}{r}} D(n,q,r,\lL)$, it suffices to show
\[D^*(n,q,r,\lL) = ((\lL/e)^{Q-1} N + o(N))^{\lL Q^{-1} \tbinom{n}{r}}.\]

\nib{Proof of Theorem \ref{count}.}
For the lower bound we start by setting aside a random subgraph $R$ of $K^r_n$ in which each edge 
is chosen with probability $n^{-a/Q}$, where we can apply Theorem \ref{hyp-decomp} 
with $a$ and $\ell$, and we suppose without loss of generality that $Q \ll h \ll \ell \ll 1/a$
(i.e.\ parameters are chosen from left to right to satisfy various inequalities below).

Next we will consider the random greedy matching 
process in the following auxiliary hypergraph $A$. We let $V(A)$ consist of $\lL$ copies of each edge 
of $K^r_n \sm R$ and $\lL-1$ copies of each edge of $R$. We let $E(A)$ consist of all 
$\tbinom{q}{r}$-sets in $V(A)$ that are edge-sets of a copy of $K^r_q$ in $K^r_n$. 

In the random greedy matching process, we start with $A$, and at each step we select a uniformly 
random edge $e$ of the current hypergraph, then delete all vertices of $e$ 
(and all incident edges) to obtain the hypergraph for the next step. 
We stop the process when fewer than $n^{r-3\ell a}$ vertices of $A$ remain and 
let $L$ denote the multigraph in $K^r_n$ consisting of the remaining vertices of $A$.
Similarly to the previous section, one can adapt the analysis of Bennett and Bohman \cite{BB} 
to show that whp (i) when $p\lL\tbinom{n}{r}$ edges remain the number of choices for 
the next edge of the process is $(1 \pm n^{-1/Q}) (p\lL)^Q \tbinom{n}{q}$, and 
(ii) $|L(e)|<2n^{1-3\ell a}$ for any $(r-1)$-set $e \sub [n]$.

Next we apply a random greedy algorithm to sequentially cover each edge of $L$ by a copy of $K^r_q$
in which all other edges are in $R$; as usual, at each step we make a uniformly random choice
subject to not using any previously covered edge. The proof of Lemma \ref{cover-leave} generalises
to show that whp the algorithm does not abort, and writing $S$ for the subgraph of $R$ covered by
the algorithm, whp $|S(e)| < n^{1-2\ell a}$ for any $(r-1)$-set $e \sub [n]$. By Chernoff bounds
whp $R$ is $(n^{-2\ell a},h)$-typical. Also whp $|R|>\tfrac{1}{2}n^{r-a/Q}$, 
so $R' := R\sm S$ is $(c_0 n^{-\ell a},h)$-typical with $|R'| > n^{r-a}$. 
Furthermore, $R'$ was obtained from $\lL K^r_n$ by deleting edge-disjoint copies of $K^r_q$,
and $\lL K^r_n$ is $K^r_q$-divisible by assumption, so $R'$ is $K^r_q$-divisible.
Therefore $R'$ has a $K^r_q$-decomposition by Theorem \ref{hyp-decomp}.
Combining this with the previous choices of $K^r_q$'s we have constructed an edge-labelled design 
with parameters $(n,q,r,\lL)$ (the vertices of $A$ specify the edge-labelling).

Now we may calculate similarly to the proof of Theorem \ref{wilson-conj}.
Writing $m$ for the number of steps in the random greedy matching process,
and $p(i) = 1-n^{-a/Q}-iQ\tbinom{n}{r}^{-1}$ for the approximate density
at the $i$th step, the logarithm of the number of choices is 
\[ L_1 = \sum_{i=1}^m (\log (p(i)^Q\lL^Q \tbinom{n}{q}) \pm 2n^{-1/Q})
= \lL Q^{-1}\tbinom{n}{r}( \log(\lL^Q \tbinom{n}{q} ) - Q \pm n^{-a/2Q} ).\]
Also, for any fixed design, the logarithm of the number of times it is counted is at most
\[L_2 =  \sum_{i=1}^m \log(p(i) \lL Q^{-1}\tbinom{n}{r}) 
= \lL Q^{-1}\tbinom{n}{r} (\log(\lL Q^{-1}\tbinom{n}{r}) - 1 \pm n^{-a/2Q} ).\]
As $\tbinom{n}{q} Q \tbinom{n}{r}^{-1} = N + o(N)$ we deduce 
\[\log D^*(n,q,r,\lL) \ge \lL Q^{-1} \tbinom{n}{r} \log ((\lL/e)^{Q-1} N + o(N)).\]

For the upper bound in Theorem \ref{count} we apply the Entropy Method,
following the argument of Linial and Luria \cite{LL} for Steiner Triple Systems
(see their paper for motivation and exposition of the method).
We let $X$ be a uniformly random edge-labelled design with parameters $(n,q,r,\lL)$, and
consider the entropy $H(X) = -\sum_D \mb{P}(X=D)\log \mb{P}(X=D)$ (using natural logarithms).
We have $D^*(n,q,r,\lL) = e^{H(X)}$, so it suffices to estimate $H(X)$.

We consider the labelled edges of $\lL K^r_n$ in a uniformly random order:
it is convenient to select $\mu = (\mu_e) \in [0,1]^{\lL K^r_n}$
uniformly at random, and to proceed by decreasing order of $\mu_e$.
At each step, when we consider $e$, we reveal the block $X_e$ of $X$
that contains $e$ and is assigned to $e$ according to the edge-labelling.
Conditional on $\mu$, we have
\[H(X) = \sum_{e \in \lL K^r_n} H(X_e \mid (X_{e'}: \mu(e')>\mu(e))).\]

We estimate $H(X_e \mid (X_{e'}: \mu(e')>\mu(e))) \le \log N^\mu_e$, where $N^\mu_e$ 
is the size of the support of the random variable $X_e \mid (X_{e'}: \mu(e')>\mu(e))$,
i.e.\ $N^\mu_e=1$ if $e$ is a labelled edge of $X_{e'}$ for some $e'$ that precedes $e$,
otherwise $N^\mu_e$ is the number of choices of a labelled $q$-set $f$ containing $e$
(i.e.\ we fix labellings of the other $Q-1$ edges in $f$) such that for each such
labelled edge $e'$ in $f$, no labelled edge of the block $X_{e'}$ precedes $e$.

Next we condition on $X$, fix $e$, and write $F_e$ for the event that $\mu(e') \le \mu(e)$ 
for all $e' \in X_e$. We estimate $\mb{E} \log N^\mu_e$, where the expectation
is with respect to $\mu$, and we suppress the $X$-conditioning in our notation.
We have \[\mb{E} \log N^\mu_e = \mb{E} (\mb{E} [ \log N^\mu_e \mid \mu_e])
= \mb{E} (\mu_e^{Q-1} \mb{E} [ \log N^\mu_e \mid \mu_e,F_e])\]
and by Jensen's inequality
$\mb{E} [ \log N^\mu_e \mid \mu_e,F_e] \le \log \mb{E} [ N^\mu_e \mid \mu_e,F_e]$.

Now we write $\mb{E} [ N^\mu_e \mid \mu_e,F_e] = 1 + \sum_f \mb{P}[E_f \mid \mu_e,F_e]$,
where the sum is over all labelled $q$-sets $f \ne X_e$ containing $e$,
and $E_f$ is the event that for each labelled edge $e'$ in $f$, 
no labelled edge of the block $X_{e'}$ precedes $e$.
Note that there are only $O(N/n)$ such $f$ with $|f \cap f'|>r$ for some block $f'$ of $X$.
For any other such $f$ we have $\mb{P}[E_f \mid \mu_e,F_e] = \mu_e^{Q(Q-1)}$.
We deduce $\mb{E} [ N^\mu_e \mid \mu_e,F_e] = \mu_e^{Q(Q-1)} \lL^{Q-1} N + O(N/n)$.

Finally, \begin{align*} & \log D^*(n,q,r,\lL) = H(X) 
\le \sum_{e \in \lL K^r_n} \mb{E} \mu_e^{Q-1} \log \mb{E} [ N^\mu_e \mid \mu_e,F_e] \\
& = \lL \tbinom{n}{r} \int_0^1 t^{Q-1} \log (t^{Q(Q-1)} \lL^{Q-1} N) + O(1/n) \ dt. \end{align*}
For any $A,B,C > 0$ we have $\int_0^1 t^{A-1} \log (Ct^B)\ dt = A^{-1}\log C - A^{-2}B$.
Setting $A=Q$, $B=Q(Q-1)$, $C=\lL^{Q-1} N$ we deduce
\begin{equation}
\log D^*(n,q,r,\lL) \le \lL Q^{-1} \tbinom{n}{r} \log ((\lL/e)^{Q-1} N + o(N)). \tag*{$\Box$}
\end{equation} 

\section{Concluding remarks}

Although we have proved (and generalised) Wilson's conjecture,
one may still ask for more precise estimates (even an asymptotic formula)
for the number of Steiner Triple Systems, and more generally designs.
Such results have been obtained by Kuperberg, Lovett and Peled \cite{KLP},
using very different methods to ours, 
but only for designs within a certain range of parameters.
One open case of particular interest (recently drawn to my attention by Ron Peled)
is the problem of estimating the number $G(n,d)$ of $d$-regular graphs on $n$ vertices. 
These may be viewed as designs with parameters $(n,2,1,d)$, 
for which our methods give $G(n,d) = d!^{-n} (dn/e + o(dn))^{dn/2}$.
Much more precise results have been obtained by McKay and Wormald,
including asymptotic enumeration for $d=\oO(n/\log n)$ (see \cite{MW1})
and $d=o(\sqrt{n})$ (see \cite{MW2}); their conjecture in \cite{MW1}
regarding a general asymptotic formula remains open.

\end{document}